\documentclass[12pt,a4paper]{article}
\usepackage{amsmath}
\usepackage{amssymb}
\usepackage{t1enc}
\usepackage[latin1]{inputenc}
\usepackage[english]{babel}
\usepackage{amsfonts}
\usepackage{latexsym}
\usepackage{bm}
\newtheorem{theorem}{Theorem}[section]

\newtheorem{lemma}[theorem]{Lemma}
\newtheorem{cor}[theorem]{Corollary}

\newtheorem{rem}[theorem]{Remark}

\begin{document}
\title{The maximum product of weights of cross-intersecting families} 

\author{Peter Borg\\[5mm]
Department of Mathematics, University of Malta\\
\texttt{peter.borg@um.edu.mt}}
\date{} 
\maketitle

\begin{abstract}
{\footnotesize Two families $\mathcal{A}$ and $\mathcal{B}$ of sets are said to be \emph{cross-$t$-intersecting} if each set in $\mathcal{A}$ intersects each set in $\mathcal{B}$ in at least $t$ elements. An active problem in extremal set theory is to determine the maximum product of sizes of cross-$t$-intersecting subfamilies of a given family. We prove a cross-$t$-intersection theorem for \emph{weighted} subsets of a set by means of a new subfamily alteration method, and use the result to provide solutions for three natural families. For $r\in[n]=\{1,2,\dots,n\}$, let ${[n]\choose r}$ be the family of $r$-element subsets of $[n]$, and let ${[n]\choose\leq r}$ be the family of subsets of $[n]$ that have at most $r$ elements. Let $\mathcal{F}_{n,r,t}$ be the family of sets in ${[n]\choose\leq r}$ that contain $[t]$. We show that if $g:{[m]\choose\leq r}\rightarrow\mathbb{R}^+$ and $h:{[n]\choose\leq s}\rightarrow\mathbb{R}^+$ are functions that obey certain conditions, $\mathcal{A}\subseteq{[m]\choose\leq r}$, $\mathcal{B}\subseteq{[n]\choose\leq s}$, and $\mathcal{A}$ and $\mathcal{B}$ are cross-$t$-intersecting, then \[\sum_{A\in\mathcal{A}}g(A)\sum_{B\in\mathcal{B}}h(B)\leq\sum_{C\in\mathcal{F}_{m,r,t}}g(C)\sum_{D\in\mathcal{F}_{n,s,t}}h(D),\] and equality holds if $\mathcal{A}=\mathcal{F}_{m,r,t}$ and $\mathcal{B}=\mathcal{F}_{n,s,t}$. We prove this in a more general setting and characterise the cases of equality. We use the result to show that the maximum product of sizes of two cross-$t$-intersecting families $\mathcal{A}\subseteq{[m]\choose r}$ and $\mathcal{B}\subseteq{[n]\choose s}$ is ${m-t\choose r-t}{n-t\choose s-t}$ for $\min\{m,n\}\geq n_0(r,s,t)$, where $n_0(r,s,t)$ is close to best possible. We obtain analogous results for families of integer sequences and for families of multisets. The results yield generalisations for $k\geq2$ cross-$t$-intersecting families, and Erdos--Ko--Rado-type results.} 
\end{abstract}

\section{Introduction} \label{Def}

Unless otherwise stated, we shall use small letters such as $x$ to
denote elements of a set or non-negative integers or functions,
capital letters such as $X$ to denote sets, and calligraphic
letters such as $\mathcal{F}$ to denote \emph{families}
(that is, sets whose elements are sets themselves). The set $\{1, 2, \dots\}$ of all positive integers is denoted by $\mathbb{N}$. 
For any $m,n \in \mathbb{N}$, the set $\{i \in \mathbb{N} \colon m \leq i \leq n\}$ is denoted by $[m,n]$. We abbreviate $[1,n]$ to $[n]$. It is to be assumed that arbitrary sets and families are \emph{finite}. We call a set $A$ an \emph{$r$-element set}, or simply an \emph{$r$-set}, if its size $|A|$ is $r$. For a set $X$, $2^X$ denotes the \emph{power set of $X$} (that is, the family of all subsets of $X$), ${X \choose r}$ denotes the family of all $r$-element subsets of $X$, and ${X \choose \leq r}$ denotes the family of all subsets of $X$ of size at most $r$. 
For a family $\mathcal{F}$ and a set $T$, we denote the family $\{F \in \mathcal{F} \colon T \subseteq F\}$ by $\mathcal{F}(T)$. 


We say that a set $A$ \emph{$t$-intersects} a set $B$ if $A$ and $B$ contain at least $t$ common elements. A family $\mathcal{A}$ of sets is said to be \emph{$t$-intersecting} if every two sets in $\mathcal{A}$ $t$-intersect. A $1$-intersecting family is also simply called an \emph{intersecting family}. If $T$ is a $t$-element subset of at least one set in a family $\mathcal{F}$, then we call the family of all the sets in $\mathcal{F}$ that contain $T$ the \emph{$t$-star of $\mathcal{F}$}. 
A $t$-star of a family is the simplest example of a $t$-intersecting subfamily.

One of the most popular endeavours in extremal set theory is that of determining the size of a largest $t$-intersecting subfamily of a given family $\mathcal{F}$. This took off with \cite{EKR}, which features the classical result, known as the Erd\H os--Ko--Rado (EKR) Theorem, that says that if $1 \leq r \leq n/2$, then the size of a largest intersecting subfamily $\mathcal{A}$ of ${[n] \choose r}$ is the size ${n-1 \choose r-1}$ of every $1$-star of ${[n] \choose r}$. If $r < n/2$, then, by the Hilton-Milner Theorem \cite{HM}, $\mathcal{A}$ attains the bound if and only if $\mathcal{A}$ is a star of ${[n] \choose r}$. If $n/2 < r \leq n$, then ${[n] \choose r}$ itself is intersecting. There are various proofs of the EKR Theorem (see \cite{Kat,HM,K,D}), two of which are particularly short and beautiful: Katona's \cite{K}, introducing the elegant cycle method, and Daykin's \cite{D}, using the fundamental Kruskal--Katona Theorem \cite{Kr,Ka}. A sequence of results \cite{EKR,F_t1,W,AK1} culminated in the solution of the problem for $t$-intersecting subfamilies of ${[n] \choose r}$; the solution particularly tells us that the size of a largest $t$-intersecting subfamily of ${[n] \choose r}$ is the size ${n-t \choose r-t}$ of a $t$-star of ${[n] \choose r}$ if and only if $n \geq (t+1)(r-t+1)$. 
The $t$-intersection problem for $2^{[n]}$ was solved by Katona \cite{Kat}. These are among the most prominent results in extremal set theory. 
The EKR Theorem inspired a wealth of results, including generalisations (see \cite{T,Borg9}), that establish how large a system of sets can be under certain intersection conditions; see \cite{DF,F,F2,Borg7,HST,HT}.

Two families $\mathcal{A}$ and $\mathcal{B}$ are said to be \emph{cross-$t$-intersecting} if each set in $\mathcal{A}$ $t$-intersects each set in $\mathcal{B}$. More generally, $k$ families $\mathcal{A}_1, \dots, \mathcal{A}_k$ (not necessarily distinct or non-empty) are said to be \emph{cross-$t$-intersecting} if for every $i$ and $j$ in $[k]$ with $i \neq j$, each set in $\mathcal{A}_i$ $t$-intersects each set in $\mathcal{A}_j$. Cross-$1$-intersecting families are also simply called \emph{cross-intersecting families}.

For $t$-intersecting subfamilies of a given family $\mathcal{F}$,
the natural question to ask is how large they can be. For
cross-$t$-intersecting families, two natural parameters arise: the
sum and the product of sizes of the cross-$t$-intersecting
families (note that the product of sizes of $k$ families
$\mathcal{A}_1, \dots, \mathcal{A}_k$ is the number of $k$-tuples
$(A_1, \dots, A_k)$ such that $A_i \in \mathcal{A}_i$ for each $i
\in [k]$). It is therefore natural to consider the problem of
maximising the sum or the product of sizes of $k$ cross-$t$-intersecting subfamilies $\mathcal{A}_1, \dots, \mathcal{A}_k$ of a given family $\mathcal{F}$. The paper \cite{Borg8} analyses this problem in general, particularly showing that for $k$ sufficiently large, both the sum and the product are maxima if $\mathcal{A}_1 = \dots = \mathcal{A}_k = \mathcal{L}$ for some largest $t$-intersecting subfamily $\mathcal{L}$ of $\mathcal{F}$. Therefore, this problem incorporates the $t$-intersection problem. Solutions have been obtained for various families (see \cite{Borg8}), including ${[n] \choose r}$ \cite{H,Pyber,MT,Bey3,Borg4,Tok1,WZ,Tok2,FLST}, 
$2^{[n]}$ \cite{MT2,Borg8}, ${[n] \choose \leq r}$ \cite{BorgBLMS}, and families of integer sequences \cite{Moon,Borg,BL2,WZ,Zhang,Tok3,FLST,PT}. 
Most of these results tell us that for the family $\mathcal{F}$ under consideration and for certain values of $k$, the sum or the product is maximum when $\mathcal{A}_1 = \dots = \mathcal{A}_k = \mathcal{L}$ for some largest $t$-star $\mathcal{L}$ of $\mathcal{F}$. In such a case, $\mathcal{L}$ is a largest $t$-intersecting subfamily of $\mathcal{F}$. 

\begin{rem} \label{remark1} \emph{In general, if $\mathcal{L} \subseteq \mathcal{F}$, $k \geq 2$, and the sum or the product is maximum when $\mathcal{A}_1 = \dots = \mathcal{A}_k = \mathcal{L}$, then $\mathcal{L}$ is a largest $t$-intersecting subfamily of $\mathcal{F}$. Indeed, the cross-$t$-intersection condition implies that every two sets $A$ and $B$ in $\mathcal{L}$ $t$-intersect (as $A \in \mathcal{A}_1$ and $B \in \mathcal{A}_2$), and by taking an arbitrary $t$-intersecting subfamily $\mathcal{A}$ of $\mathcal{F}$ and setting $\mathcal{B}_1 = \dots = \mathcal{B}_k = \mathcal{A}$, we obtain that $\mathcal{B}_1, \dots, \mathcal{B}_k$ are cross-$t$-intersecting, and hence $|\mathcal{A}| \leq |\mathcal{L}|$ since $k|\mathcal{A}| = \sum_{i=1}^k |\mathcal{B}_i| \leq \sum_{i=1}^k |\mathcal{A}_i| = k|\mathcal{L}|$ or $|\mathcal{A}|^k = \prod_{i=1}^k |\mathcal{B}_i| \leq \prod_{i=1}^k |\mathcal{A}_i| = |\mathcal{L}|^k$.}
\end{rem}

Wang and Zhang \cite{WZ} solved the maximum sum problem for an important class of families that includes ${[n] \choose r}$ and families of integer sequences, using a striking combination of the method in \cite{Borg4,Borg3,Borg2,BL2,Borg5} and an important lemma that is found in \cite{AC,CK} and referred to as the `no-homomorphism lemma'. The solution for ${[n] \choose r}$ with $t=1$ had been obtained by Hilton \cite{H} and is the first result of this kind. 

In this paper we address the maximum product problem for ${[n] \choose r}$ and families of integer sequences. We will actually consider more general problems; one generalisation allows the cross-$t$-intersecting families to come from different families, and another one involves maximising instead the product of \emph{weights} of cross-$t$-intersecting families of subsets of a set. 
As we explain in the next section, if the product for $k = 2$ is maximum when the cross-$t$-intersecting families are certain $t$-stars, then this immediately generalises for $k \geq 2$. 

The maximum product problem for ${[n] \choose r}$ was first addressed by Pyber \cite{Pyber}, who proved that for any $r$, $s$, and $n$ such that either $r = s \leq n/2$ or $r < s$ and $n \geq 2s + r -2$, if $\mathcal{A} \subseteq {[n] \choose r}$ and $\mathcal{B} \subseteq {[n] \choose s}$ such that $\mathcal{A}$ and $\mathcal{B}$ are cross-intersecting, then $|\mathcal{A}||\mathcal{B}| \leq {n-1 \choose r-1}{n-1 \choose s-1}$. Subsequently, Matsumoto and Tokushige \cite{MT} proved this for $r \leq s \leq n/2$. 
It has been shown in \cite{Borg11} that there exists an integer $n_0(r,s,t)$ such that for $t \leq r \leq s$ and $n \geq n_0(r,s,t)$, if $\mathcal{A} \subseteq {[n] \choose r}$, $\mathcal{B} \subseteq {[n] \choose s}$, and $\mathcal{A}$ and $\mathcal{B}$ are cross-$t$-intersecting, then $|\mathcal{A}||\mathcal{B}| \leq {n-t \choose r-t}{n-t \choose s-t}$. The value of $n_0(r,s,t)$ given in \cite{Borg11} is far from best possible. The special case $r = s$ is treated in \cite{Tok1,Tok2,FLST}, which establish values of $n_0(r,r,t)$ that are close to the conjectured smallest value of $(t+1)(r-t+1)$, and which use algebraic methods and Frankl's random walk method \cite{F_t1}; in particular, $n_0(r,r,t) = (t+1)r$ is determined in \cite{FLST} for $t \geq 14$. 
Using purely combinatorial arguments, we solve the problem for $n \geq (t+u+2)(s-t)+r-1$, where $u$ can be any non-negative real number satisfying $u > \frac{6-t}{3}$; thus, we can take $n_0(r,s,t) = (t+2)(s-t)+r-1$ for $t \geq 7$, and $n_0(r,s,t) < (t+4)(s-t)+r-1$ for $1 \leq t \leq 6$. We actually prove the following more general result in Section~\ref{nrssection}.

\begin{theorem}\label{nrs} If $1 \leq t \leq r \leq s$, $u$ is a non-negative real number such that $u > \frac{6-t}{3}$, $\min\{m,n\} \geq (t+u+2)(s-t)+r-1$, $\mathcal{A} \subseteq {[m] \choose r}$, $\mathcal{B} \subseteq {[n] \choose s}$, and $\mathcal{A}$ and $\mathcal{B}$ are cross-$t$-intersecting, then
\[|\mathcal{A}||\mathcal{B}| \leq {m-t \choose r-t}{n-t \choose s-t}.\]
Moreover, if $u > 0$, then the bound is attained if and only if $\mathcal{A} = \left\{ A \in {[m] \choose r} \colon T \subseteq A \right\}$ and $\mathcal{B} = \left\{ B \in {[n] \choose s} \colon T \subseteq B \right\}$ for some $t$-element subset $T$ of $[\min\{m,n\}]$.
\end{theorem}
%

In Section~\ref{nrssection}, we show that Theorem~\ref{nrs} is a consequence of our main result, Theorem~\ref{xintweight}, for which we need some additional definitions and notation. 

For any $i, j \in [n]$, let $\delta_{i,j} \colon 2^{[n]} \rightarrow 2^{[n]}$ be defined by
\[ \delta_{i,j}(A) = \left\{ \begin{array}{ll}
(A \backslash \{j\}) \cup \{i\} & \mbox{if $j \in A$ and $i \notin
A$};\\
A & \mbox{otherwise,}
\end{array} \right. \]
and let $\Delta_{i,j} \colon 2^{2^{[n]}} \rightarrow 2^{2^{[n]}}$ be the \emph{compression operation} defined by
\[\Delta_{i,j}(\mathcal{A}) = \{\delta_{i,j}(A) \colon A \in
\mathcal{A}\} \cup \{A \in \mathcal{A} \colon \delta_{i,j}(A) \in \mathcal{A}\}.\]
The compression operation was introduced in the seminal paper \cite{EKR}. The paper \cite{F} provides a survey on the properties and uses of compression (also called \emph{shifting}) operations in extremal set theory. All our new results make use of compression operations.

If $i < j$, then we call $\Delta_{i,j}$ a \emph{left-compression}. A family $\mathcal{F} \subseteq 2^{[n]}$ is said to be
\emph{compressed} if $\Delta_{i,j}(\mathcal{F}) = \mathcal{F}$ for every $i,j \in [n]$ with $i < j$. In other words, $\mathcal{F}$ is compressed if it is invariant under left-compressions. Note that $\mathcal{F}$ is compressed if and only if $(F \backslash \{j\}) \cup \{i\} \in \mathcal{F}$ whenever $i < j \in F \in \mathcal{F}$ and $i \in [n] \backslash F$.

A family $\mathcal{H}$ is said to be \emph{hereditary} if for each $H \in \mathcal{H}$, all the subsets of $H$ are in $\mathcal{H}$. Thus,
a family is hereditary if and only if it is a union of power sets. The family ${[n] \choose \leq r}$ (which is $2^{[n]}$ if $r = n$) is an example of a hereditary family that is compressed.

Let $\mathbb{R}^+$ denote the set of positive real numbers. With a slight abuse of notation, for any non-empty family $\mathcal{F}$, any function $w \colon \mathcal{F} \rightarrow \mathbb{R}^+$ (called a \emph{weight function}), and any $\mathcal{A} \subseteq \mathcal{F}$, we denote the sum $\sum_{A \in \mathcal{A}} w(A)$ (of \emph{weights} of
sets in $\mathcal{A}$) by $w(\mathcal{A})$. Note that if $\mathcal{A}$ is empty, then $w(\mathcal{A})$ is the \emph{empty sum}, and we will adopt the convention of taking this to be $0$.

In Section~\ref{Weightedsection}, we prove the following result.

\begin{theorem}\label{xintweight} Let $1 \leq t \leq n$, $T = [t]$, and $u \in \{0\} \cup \mathbb{R}^+$ such that $u > \frac{6-t}{3}$. Let $\mathcal{G}$ and $\mathcal{H}$ be non-empty compressed hereditary subfamilies of $2^{[n]}$. For each $\mathcal{F} \in \{\mathcal{G}, \mathcal{H}\}$, let $w_{\mathcal{F}} \colon \mathcal{F} \rightarrow \mathbb{R}^+$ such that \\
(a) $w_{\mathcal{F}}(A) \geq (t+u)w_{\mathcal{F}}(B)$ for every $A, B \in \mathcal{F}$ with $A \subsetneq B$ and $|A| \geq t$, and \\
(b) $w_{\mathcal{F}}(\delta_{i,j}(C)) \geq w_{\mathcal{F}}(C)$ for every $C \in \mathcal{F}$ and every $i,j \in [n]$ with $i < j$. \\
Let $g = w_{\mathcal{G}}$ and $h = w_{\mathcal{H}}$. If $\mathcal{A} \subseteq \mathcal{G}$ and $\mathcal{B} \subseteq  \mathcal{H}$ such that $\mathcal{A}$ and $\mathcal{B}$ are cross-$t$-intersecting, then
\[g(\mathcal{A}) h(\mathcal{B}) \leq g(\mathcal{G}(T)) h(\mathcal{H}(T)).\]
Moreover, if $u > 0$ and each of $\mathcal{G}$ and $\mathcal{H}$ has a member of size at least $t$, then the bound is attained if and only if $\mathcal{A} = \mathcal{G}(T')$ and $\mathcal{B} = \mathcal{H}(T')$ for some $T' \in {[n] \choose t}$ such that $g(\mathcal{G}(T')) = g(\mathcal{G}(T))$ and $h(\mathcal{H}(T')) = h(\mathcal{H}(T))$.  
\end{theorem}
%

\begin{rem}\label{remark2} \emph{For $u > \frac{6-t}{3}$ to hold, we can always take $u = 2$, and we can take $u = 0$ for $t \geq 7$. We conjecture that the inequality $g(\mathcal{A}) h(\mathcal{B}) \leq g(\mathcal{G}(T)) h(\mathcal{H}(T))$ still holds if the condition $u > \frac{6-t}{3}$ is replaced by $u=0$. As we mentioned above, this is true for $t \geq 7$. Also, the proof of Theorem~\ref{xintweight} shows that for $t \geq 3$, the conjecture is true if it is true for $t+3 \leq n \leq t+6$ (see Remark~\ref{proofremark}). A verification of the conjecture for $t+3 \leq n \leq t+6$ could be obtained through detailed case-checking similar to that used in our proof for the special case $n \leq t+2$; however, the process would be significantly more laborious. 
The condition on $u$ cannot be relaxed further, because no real number $u < 0$ with $t+u \geq 1$ guarantees that the result holds. Indeed, if $1 \leq x = t+u < t \leq n-2$, $\mathcal{G} = \mathcal{H} = 2^{[n]}$, $g(G) = h(G) = x^{n-|G|}$ for all $G \in 2^{[n]}$, and $\mathcal{A} = \mathcal{B} = \{A \in 2^{[n]} \colon |A \cap [t+2]| \geq t+1\} = \mathcal{A}^*$, then conditions (a) and (b) of Theorem~\ref{xintweight} are satisfied, $\mathcal{A}$ and $\mathcal{B}$ are cross-$t$-intersecting, but 
\begin{align} (g(\mathcal{A})h(\mathcal{B}))^{1/2} &= g(\mathcal{A}) = \sum_{X \in {[t+2] \choose t+1}}\sum_{Y \subseteq {[t+3,n]}} g(X \cup Y) + \sum_{Y \subseteq {[t+3,n]}} g([t+2] \cup Y) \nonumber \\
&= (t+2)\sum_{j = 0}^{n-t-2} {n-t-2 \choose j} x^{n-t-1-j} + \sum_{j = 0}^{n-t-2}{n-t-2 \choose j}x^{n-t-2-j} \nonumber \\
&= (t+2)x^{n-t-1}\left(1 + x^{-1} \right)^{n-t-2} + x^{n-t-2}\left(1 + x^{-1} \right)^{n-t-2} \nonumber \\
&= x^{n-t-2}\left(1 + x^{-1} \right)^{n-t-2} (tx + 2x + 1) = (x+1)^{n-t-2}(tx + 2x + 1) \nonumber \\
&> (x+1)^{n-t-2} (x^2 + 2x + 1) \quad \mbox{(as $1 \leq x < t$)} \nonumber \\
&= x^{n-t}\left(1 + x^{-1} \right)^{n-t} = \sum_{j = 0}^{n-t}{n-t \choose j}x^{n-t-j} = \sum_{Y \subseteq {[t+1,n]}} g([t] \cup Y) \nonumber \\
&= g(\mathcal{G}(T)) = (g(\mathcal{G}(T))h(\mathcal{H}(T)))^{1/2}, \nonumber 
\end{align}
and hence $g(\mathcal{A})h(\mathcal{B}) > g(\mathcal{G}(T))h(\mathcal{H}(T))$. It has been shown in \cite{BorgBLMS} that for $t = 1$, the product of sizes of $\mathcal{A}$ and $\mathcal{B}$ is maximised by taking $\mathcal{A} = \mathcal{G}(T)$ and $\mathcal{B} = \mathcal{H}(T)$; equivalently, for the special case where $t = 1$ and $g(A) = h(A) = 1$ for all $A \in \mathcal{G} \cup \mathcal{H}$, the bound in Theorem~\ref{xintweight} also holds (that is, the conjecture is true). 
However, this is not true for $t > 1$, and hence Theorem~\ref{xintweight} does not imply that the product of sizes is maximised by taking $\mathcal{A} = \mathcal{G}(T)$ and $\mathcal{B} = \mathcal{H}(T)$. Indeed, if $\mathcal{G} = \mathcal{H} = 2^{[n]}$ and $\mathcal{A} = \mathcal{B} = \mathcal{A}^*$ as above, then $|\mathcal{A}||\mathcal{B}| > |\mathcal{G}(T)| |\mathcal{H}(T)|$ (take $x = 1$ above).} 
\end{rem}

The proof of Theorem~\ref{xintweight} contains the main observations in this paper and is based on induction, compression, a new subfamily alteration method, and double-counting. The alteration method can be regarded as the main new component and appears to have the potential of yielding other intersection results of this kind. 

The bound in \cite[Theorem~1.3]{FLST} for product measures of cross-$t$-intersecting subfamilies of $2^{[n]}$ is given by Theorem~\ref{xintweight} 
with $\mathcal{G} = \mathcal{H} = 2^{[n]}$, $t \geq 14$, $u = 0$, and $g(A) = h(A) = p^{|A|}(1-p)^{n-|A|}$ for all $A \in 2^{[n]}$, where $p \in \mathbb{R}^+$ such that $p \leq \frac{1}{t+1}$. 

The subsequent results in this section and in the next section are also consequences of Theorem~\ref{xintweight}. Our next application is a cross-$t$-intersection result for integer sequences. 
 
We will represent a sequence $a_1, \dots, a_n$ by an $n$-tuple $(a_1, \dots, a_n)$, and we say that it is of \emph{length $n$}. 
We call a sequence of positive integers a \emph{positive sequence}. We call $(a_1, \dots, a_n)$ an \emph{$r$-partial sequence} if exactly $r$ of its entries are positive integers and the rest are all zero. Thus, an $n$-partial sequence of length $n$ is positive. A sequence $(c_1, \dots, c_n)$ is said to be \emph{increasing} if $c_1 \leq \dots \leq c_n$. We call an increasing positive sequence an \emph{IP sequence}. Note that $(c_1, \dots, c_n)$ is an IP sequence if and only if $1 \leq c_1 \leq \dots \leq c_n$. 

We call $\{(x_1,y_1), \dots, (x_r,y_r)\}$ a \emph{labeled set} (following \cite{Borg}) if $x_1, \dots, x_r$ are distinct. For any IP sequence ${\bf c} = (c_1, \dots, c_n)$ and any $r \in [n]$, let $\mathcal{S}_{{\bf c},r}$ be the family of all labeled sets $\{(x_1,y_{x_1}), \dots, (x_r, y_{x_r})\}$ such that $\{x_1, \dots, x_r\} \in {[n] \choose r}$ and $y_{x_j} \in [c_{x_j}]$ for each $j \in [r]$. 
%
%
For any sets $Y_1, \dots, Y_n$, let $Y_1 \times \dots \times Y_n$ denote the \emph{Cartesian product of $Y_1, \dots, Y_n$}, that is, the set of sequences $(y_1, \dots, y_n)$ such that $y_i \in Y_i$ for each $i \in [n]$. 
Note that $\mathcal{S}_{{\bf c},n} = \{\{(1,y_1), \dots, (n,y_n)\} \colon y_i \in [c_i] \mbox{ for each } i \in [n]\}$, so $\mathcal{S}_{{\bf c},n}$ is isomorphic to $[c_1] \times \dots \times [c_n]$. Also note that $\mathcal{S}_{{\bf c},r}$ is isomorphic to the set of $r$-partial sequences $(y_1, \dots, y_n)$ such that for some $R \in {[n] \choose r}$, $y_i \in [c_i]$ for each $i \in R$ (and hence $y_j = 0$ for each $j \in [n] \backslash R$). Let $\mathcal{S}_{{\bf c},r,t} = \mathcal{S}_{{\bf c},r}([t] \times [1]) = \left\{A \in \mathcal{S}_{{\bf c},r} \colon (x,1) \in A \mbox{ for each } x \in [t]\right\}$. 

In Section~\ref{Proofmain}, we prove the following result.

\begin{theorem} \label{main} Let ${\bf c} = (c_1, \dots, c_m)$ and ${\bf d} = (d_1, \dots, d_n)$ be IP sequences. Let $r \in [m]$, $s \in [n]$, $t \in [\min\{r,s\}]$, and $u \in \{0\} \cup \mathbb{R}^+$ such that $u > \frac{6-t}{3}$. If $c_1 \geq t+u+1$, $d_1 \geq t+u+1$, $\mathcal{A} \subseteq \mathcal{S}_{{\bf c},r}$, $\mathcal{B} \subseteq \mathcal{S}_{{\bf d},s}$, and $\mathcal{A}$ and $\mathcal{B}$ are cross-$t$-intersecting, then
\[|\mathcal{A}||\mathcal{B}| \leq \Bigg{(} \sum_{I \in {[t+1,m] \choose r-t}} \prod_{i \in I} c_i \Bigg{)} \Bigg{(} \sum_{J \in {[t+1,n] \choose s-t}}\prod_{j \in J} c_j \Bigg{)}.\]
Moreover, if $u > 0$, then the bound is attained if and only if for some $T \in \mathcal{S}_{{\bf c},t} \cap \mathcal{S}_{{\bf d},t}$ with $|\mathcal{S}_{{\bf c},r}(T)| = |\mathcal{S}_{{\bf c},r,t}|$ and $|\mathcal{S}_{{\bf d},s}(T)| = |\mathcal{S}_{{\bf d},s,t}|$, $\mathcal{A} = \mathcal{S}_{{\bf c},r}(T)$ and $\mathcal{B} = \mathcal{S}_{{\bf d},s}(T)$.
\end{theorem}
%
%
%
%
Note that this result holds for $c_1 \geq t+1$ and $d_1 \geq t+1$ when $t \geq 7$, for $c_1 \geq t+2$ and $d_1 \geq t+2$ when $4 \leq t \leq 6$, and for $c_1 \geq t+3$ and $d_1 \geq t+3$ when $1 \leq t \leq 3$. We conjecture that the result holds for $c_1 \geq t+1$ and $d_1 \geq t+1$, and, as can be seen from the proof of Theorem~\ref{main}, this conjecture is true if the conjecture in Remark~\ref{remark2} is true. The result does not hold for $c_1 < t+1$. Indeed, if $r = s = m = n \geq t+2$, $c_1 = \dots = c_n = x+1 < t+1 = d_1 = \dots = d_n$, $Z = [n] \times [1]$, $Z_1 = [t+2] \times [1]$, $Z_2 = [t+3,n] \times [1]$, $\mathcal{A} = \{A \in \mathcal{S}_{{\bf c},n} \colon |A \cap Z_1| \geq t+1\}$, and $\mathcal{B} = \{B \in \mathcal{S}_{{\bf d},n} \colon |B \cap Z_1| \geq t+1\}$, then $\mathcal{A}$ and $\mathcal{B}$ are cross-$t$-intersecting, 
\begin{align} |\mathcal{A}| &= \left|\bigcup_{X \in {Z_1 \choose t+1} \cup \{Z_1\}}\bigcup_{j=0}^{|Z_2|}\bigcup_{Y \in {Z_2 \choose j}} \{A \in \mathcal{A} \colon A \cap Z = X \cup Y\} \right| \nonumber \\
&= (t+2)\sum_{j = 0}^{n-t-2} {n-t-2 \choose j} x^{n-t-1-j} + \sum_{j = 0}^{n-t-2}{n-t-2 \choose j}x^{n-t-2-j} \nonumber \\
&= \left(x + 1 \right)^{n-t-2} (tx + 2x + 1) \quad \mbox{(as in Remark~\ref{remark2})} \nonumber \\
&> \left(x + 1 \right)^{n-t-2} (x^2 + 2x + 1) = (x+1)^{n-t} = |\mathcal{S}_{{\bf c},n,t}|, \nonumber 
\end{align}  
$|\mathcal{B}| = \left(t + 1 \right)^{n-t-2} (t^2 + 2t + 1) = (t+1)^{n-t} = |\mathcal{S}_{{\bf d},n,t}|$ (by a calculation similar to that for $|\mathcal{A}|$), and hence $|\mathcal{A}||\mathcal{B}| > |\mathcal{S}_{{\bf c},r,t}||\mathcal{S}_{{\bf d},s,t}|$.

Solutions for the special case where ${\bf c} = {\bf d}$ and $r = s = n$ already exist. The solution for $t+2 \leq c_1 = c_n$ was first obtained by Moon \cite{Moon}. Inspired by \cite{Zhang}, Pach and Tardos \cite{PT} recently generalised Moon's result to include the cases $t+2 \leq c_1 \leq c_n$ and $8 \leq t+1 \leq c_1 \leq c_n$. Another proof for $15 \leq t+1 \leq c_1 = c_n$ is given in \cite{FLST}. 
%
 
Our last application of Theorem~\ref{xintweight} in this section is a cross-$t$-intersection result for multisets. 

A \emph{multiset} is a collection $A$ of objects such that each object possibly appears more than once in $A$. Thus the difference between a multiset and a set is that a multiset may have repetitions of its elements. We can uniquely represent a multiset $A$ of positive integers by an IP sequence $(a_1, \dots, a_r)$, where $a_1, \dots, a_r$ form $A$. Thus we will take multisets to be IP sequences. For $A = (a_1, \dots, a_r)$, the \emph{support of $A$} is the set $\{a_1, \dots, a_r\}$ and will be denoted by ${\rm S}_A$. For any $n, r \in \mathbb{N}$, let $M_{n,r}$ denote the set of all multisets $(a_1, \dots, a_r)$ such that $a_1, \dots, a_r \in [n]$; thus $M_{n,r} = \{(a_1, \dots, a_r) \colon a_1 \leq \dots \leq a_r, \, a_1, \dots, a_r \subseteq [n]\}$. An elementary counting result is that \[|M_{n,r}| = {n+r-1 \choose r}.\]

With a slight abuse of terminology, we say that a multiset $A$ \emph{$t$-intersects} a multiset $B$ if and $A$ and $B$ have at least $t$ distinct common elements, that is, if ${\rm S}_A$ $t$-intersects ${\rm S}_B$. A set $\mathcal{A}$ of multisets is said to be \emph{$t$-intersecting} if every two multisets in $\mathcal{A}$ $t$-intersect, and $k$ sets $\mathcal{A}_1, \dots, \mathcal{A}_k$ of multisets are said to be \emph{cross-$t$-intersecting} if for every $i,j \in [k]$ with $i \neq j$, each multiset in $\mathcal{A}_i$ $t$-intersects each multiset in $\mathcal{A}_j$.

In Section~\ref{multisection}, we prove the following result.

\begin{theorem}\label{multithm} If $1 \leq t \leq r \leq s$, $u \in \{0\} \cup \mathbb{R}^+$ such that $u > \frac{6-t}{3}$, $\min\{m,n\} \geq (t+u+1)(s-t)+r-t$, $\mathcal{A} \subseteq M_{m,r}$, $\mathcal{B} \subseteq M_{n,s}$, and $\mathcal{A}$ and $\mathcal{B}$ are cross-$t$-intersecting, then
\[|\mathcal{A}||\mathcal{B}| \leq {m+r-t-1 \choose r-t}{n+s-t-1 \choose s-t}.\]
Moreover, if $u > 0$, then the bound is attained if and only if $\mathcal{A} = \left\{ A \in M_{m,r} \colon T \subseteq {\rm S}_A \right\}$ and $\mathcal{B} = \left\{ B \in M_{n,s} \colon T \subseteq {\rm S}_B \right\}$ for some $t$-element subset $T$ of $[\min\{m,n\}]$.
\end{theorem}
The condition $\min\{m,n\} \geq (t+u+1)(s-t)+r-t$ is close to being sharp, as is evident from the fact that if $r = s$, $m = n < t(r-t)+2$, and $\mathcal{A} = \mathcal{B} = \{A \in M_{n,r} \colon |{\rm S}_A \cap [t+2]| \geq t+1\}$, then $\mathcal{A}$ and $\mathcal{B}$ are cross-$t$-intersecting,

\begin{align} 
|\mathcal{A}| &= \sum_{X \in {[t+2] \choose t+1} \cup \{[t+2]\}} |\{A \in M_{n,r} \colon {\rm S}_A \cap [t+2] = X\}| \nonumber \\
&= \sum_{X \in {[t+2] \choose t+1} \cup \{[t+2]\}} |\{(a_1, \dots, a_{r-|X|}) \colon a_1 \leq  \dots \leq a_{r-|X|}, \, a_1, \dots, a_{r-|X|} \in X \cup [t+3,n]\}| \nonumber \\
&= \sum_{X \in {[t+2] \choose t+1} \cup \{[t+2]\}} |M_{|X|+n-t-2,r-|X|}| \nonumber \\
&= \sum_{X \in {[t+2] \choose t+1} \cup \{[t+2]\}} {n+r-t-3 \choose r-|X|} = (t+2){n+r-t-3 \choose r-t-1} + {n+r-t-3 \choose r-t-2} \nonumber \\
&= \frac{{n+r-t-1 \choose r-t}}{(n+r-t-1)(n+r-t-2)}\left( (t+2)(r-t)(n-1) + (r-t)(r-t-1) \right) \nonumber \\
&> \frac{{n+r-t-1 \choose r-t}}{((t+1)(r-t)+1)((t+1)(r-t))}\left( (t+2)(r-t)(t(r-t)+1) + (r-t)(r-t-1) \right) \nonumber \\
&= {n+r-t-1 \choose r-t}, \nonumber
\end{align}
and hence $|\mathcal{A}||\mathcal{B}| > {m+r-t-1 \choose r-t}^2 = {m+r-t-1 \choose r-t}{n+s-t-1 \choose s-t}$. 

EKR-type results for multisets have been obtained in \cite{MP,FGV}. To the best of the author's knowledge, Theorem~\ref{multithm} is the first cross-$t$-intersection result for multisets.
 
In the next section, we show that the above results generalise for $k \geq 2$ families and yield EKR-type results. Section~\ref{Compsection} provides basic compression results used in our proofs. 
Sections~\ref{Weightedsection}--\ref{multisection} are dedicated to the proofs of Theorems~\ref{xintweight}, \ref{nrs}, \ref{main}, and \ref{multithm}, respectively. 

\section{Multiple cross-$t$-intersecting families and $t$-intersecting families}

Theorem~\ref{nrs} generalises as follows.

\begin{theorem} \label{maincor} Let $k \geq 2$, $t \leq r_1 \leq \dots \leq r_k$, $u \in \{0\} \cup \mathbb{R}^+$ such that $u > \frac{6-t}{3}$, and $\min\{n_1, \dots, n_k\} \geq (t+u+2)(r_k-t)+r_{k-1}-1$. If $\mathcal{A}_1 \subseteq {[n_1] \choose r_1}, \dots, \mathcal{A}_k \subseteq {[n_k] \choose r_k}$, and $\mathcal{A}_1, \dots, \mathcal{A}_k$ are cross-$t$-intersecting, then
\[\prod_{i=1}^k |\mathcal{A}_i| \leq \prod_{i=1}^k {n_i-t \choose r_i - t}.\]
Moreover, if $u > 0$, then the bound is attained if and only if for some $t$-element subset $T$ of $[\min\{n_1, \dots, n_k\}]$, $\mathcal{A}_i = \{A \in {[n_i] \choose r_i} \colon T \subseteq A\}$ for each $i \in [k]$.
\end{theorem}
%
The line of argument in the proof of \cite[Theorem~1.2]{Borg11} yields the result above together with a similar generalisation of Theorem~\ref{multithm} and the following generalisations of Theorem~\ref{xintweight} and Theorem~\ref{main}. 

\begin{theorem}\label{xintweightgen} If $t, u$, and $T$ are as in Theorem~\ref{xintweight}, $\mathcal{H}_1, \dots, \mathcal{H}_k$ are non-empty compressed hereditary subfamilies of $2^{[n]}$, $w_{\mathcal{F}} \colon \mathcal{F} \rightarrow \mathbb{R}^+$ is a function satisfying (a) and (b) (of Theorem~\ref{xintweight}) for each $\mathcal{F} \in \{\mathcal{H}_1, \dots, \mathcal{H}_k\}$, $\mathcal{A}_i \subseteq \mathcal{H}_i$ for each $i \in [k]$, and $\mathcal{A}_1, \dots, \mathcal{A}_k$ are cross-$t$-intersecting, then 
\[\prod_{i=1}^k w_{\mathcal{H}_i}(\mathcal{A}_i) \leq \prod_{i=1}^k w_{\mathcal{H}_i}(\mathcal{H}_i(T)).\]
Moreover, if $u > 0$ and each of $\mathcal{H}_1, \dots, \mathcal{H}_k$ has a member of size at least $t$, then the bound is attained if and only if for some $T' \in {[n] \choose t}$ such that $w_{\mathcal{H}_i}(\mathcal{H}_i(T')) = w_{\mathcal{H}_i}(\mathcal{H}_i(T))$ for each $i \in [k]$, $\mathcal{A}_i = \mathcal{H}_i(T')$ for each $i \in [k]$.
\end{theorem} 

\begin{theorem} Let ${\bf c}_1 = (c_{1,1}, \dots, c_{1,n_1}), \dots, {\bf c}_k = (c_{k,1}, \dots, c_{k,n_k})$ be IP sequences. Let $r_1 \in [n_1], \dots, r_k \in [n_k]$, $t \in [\min\{r_1, \dots, r_k\}]$, and $u \in \{0\} \cup \mathbb{R}^+$ such that $u > \frac{6-t}{3}$. If $c_{1,1} \geq t+u+1, \dots, c_{k,1} \geq t+u+1$, $\mathcal{A}_1 \subseteq \mathcal{S}_{{\bf c}_1,r_1}, \dots, \mathcal{A}_k \subseteq \mathcal{S}_{{\bf c}_k,r_k}$, and $\mathcal{A}_1, \dots, \mathcal{A}_k$ are cross-$t$-intersecting, then
\[\prod_{i=1}^k |\mathcal{A}_i| \leq \prod_{i=1}^k \Bigg{(} \sum_{I \in {[t+1,n_i] \choose r_i-t}} \prod_{j \in I} c_{i,j} \Bigg{)}.\]
Moreover, if $u > 0$, then the bound is attained if and only if for some $T \in \bigcap_{i=1}^k \mathcal{S}_{{\bf c}_i,t}$ with $|\mathcal{S}_{{\bf c}_i,r_i}(T)| = |\mathcal{S}_{{\bf c}_i,r_i,t}|$ for each $i \in [k]$, $\mathcal{A}_i = \mathcal{S}_{{\bf c}_i,r_i}(T)$ for each $i \in [k]$.
\end{theorem}
We simply observe that $\left( \prod_{i=1}^k a_i \right)^{k-1} = \prod_{i=1}^k \prod_{j \in [k] \backslash [i]} a_ia_j$ (see also \cite[Lemma~5.2]{Borg8} with $p = 2$) and that if $\mathcal{A}_1, \dots, \mathcal{A}_k$ are cross-$t$-intersecting, then any $\mathcal{A}_i$ and $\mathcal{A}_j$ with $i \neq j$ are cross-$t$-intersecting.  Thus, if, for example, $\mathcal{A}_1, \dots, \mathcal{A}_k$ are as in Theorem~\ref{xintweightgen}, $a_i = w_{\mathcal{H}_i}(\mathcal{A}_i)$ for each $i \in [k]$, and $b_i = w_{\mathcal{H}_i}(\mathcal{H}_i(T))$ for each $i \in [k]$, then Theorem~\ref{xintweight} gives us $\prod_{i=1}^k \prod_{j \in [k] \backslash [i]} a_ia_j \leq \prod_{i=1}^k \prod_{j \in [k] \backslash [i]} b_ib_j$, and hence $\left( \prod_{i=1}^k a_i \right)^{k-1} \leq \left( \prod_{i=1}^k b_i \right)^{k-1}$ (giving $\prod_{i=1}^k a_i \leq \prod_{i=1}^k b_i$, as required). 

As in Remark~\ref{remark1}, Theorem~\ref{xintweight} immediately implies an EKR-type version for a family $\mathcal{H}$ as in Theorem~\ref{xintweight}. By taking $\mathcal{G} = \mathcal{H}$ in Theorem~\ref{xintweight} and applying an argument similar to the one in Remark~\ref{remark1}, we obtain the following new result.

\begin{theorem} Let $t, u, T, \mathcal{H}$, and $h$ be as in Theorem~\ref{xintweight}. If $\mathcal{A}$ is a $t$-intersecting subfamily of $\mathcal{H}$, then 
\[h(\mathcal{A}) \leq h(\mathcal{H}(T)).\] 
Moreover, if $u > 0$ and $\mathcal{H}$ has a member of size at least $t$, then the bound is attained if and only if $\mathcal{A} = \mathcal{H}(T')$ for some $t$-set $T'$ such that $h(\mathcal{H}(T')) = h(\mathcal{H}(T))$.
\end{theorem} 

By taking ${\bf c} = {\bf d}$ in Theorem~\ref{main} and applying the argument in Remark~\ref{remark1}, we obtain the following EKR-type result. 

\begin{theorem} \label{mainEKR} If $1 \leq t \leq r \leq n$, $u \in \{0\} \cup \mathbb{R}^+$ such that $u > \frac{6-t}{3}$, ${\bf c} = (c_1, \dots, c_n)$ is an IP sequence, $c_1 \geq t+u+1$, and $\mathcal{A}$ is a $t$-intersecting subfamily of $\mathcal{S}_{{\bf c},r}$, then
\[|\mathcal{A}| \leq \Bigg{(} \sum_{I \in {[t+1,n] \choose r-t}} \prod_{i \in I} c_i \Bigg{)}.\]
Moreover, if $u > 0$, then the bound is attained if and only if $\mathcal{A} = \mathcal{S}_{{\bf c},r}(T)$ for some $T \in \mathcal{S}_{{\bf c},t}$ with $|\mathcal{S}_{{\bf c},r}(T)| = |\mathcal{S}_{{\bf c},r,t}|$.
\end{theorem}
The EKR problem for $\mathcal{S}_{{\bf c},r}$ attracted much attention and has been dealt with extensively (see, for example, \cite{Borg7}). In particular, for $c_1 = c_n$, it was solved for $r = n$ in \cite{AK2,FT2}, and for $n \geq \left \lfloor \frac{(r-t+c_1)(t+1)}{c_1} \right \rfloor$ in \cite{Bey1}. 
Similarly to Theorem~\ref{main}, Theorem~\ref{mainEKR} does not hold for $c_1 < t+1$.

By taking $m = n$ and $r = s$ in Theorem~\ref{multithm}, and applying the argument in Remark~\ref{remark1}, we obtain the following EKR-type result. 

\begin{theorem}\label{multigen} If $1 \leq t \leq r$, $u \in \{0\} \cup \mathbb{R}^+$ such that $u > \frac{6-t}{3}$, $n \geq (t+u+2)(r-t)$, $\mathcal{A} \subseteq M_{n,r}$, and $\mathcal{A}$ is $t$-intersecting, then
\[|\mathcal{A}| \leq {n+r-t-1 \choose r-t}.\]
Moreover, if $u > 0$, then the bound is attained if and only if $\mathcal{A} = \left\{ A \in M_{n,r} \colon T \subseteq {\rm S}_A \right\}$ for some $T \in {[n] \choose t}$.
\end{theorem}
The condition $n \geq (t+u+2)(r-t)$ is close to being sharp. Indeed, as shown in Section~\ref{Def}, if $n < t(r-t) + 2$ and $\mathcal{A} = \{A \in M_{n,r} \colon |{\rm S}_A \cap [t+2]| \geq t+1\}$, then $|\mathcal{A}| > {n+r-t-1 \choose r-t}$.

The EKR problem for $M_{n,r}$ and $t=1$ is solved in \cite{MP}. Generalising this result, F\"{u}redi, Gerbner, and Vizer \cite{FGV} solved the EKR problem of maximising the size of a largest subset $\mathcal{A}$ of $M_{n,r}$ such that for every $(a_1, \dots, a_r), (b_1, \dots, b_r) \in \mathcal{A}$, there exist $t$ distinct elements $i_1, \dots, i_t$ of $[r]$ and $t$ distinct elements $j_1, \dots, j_t$ of $[r]$ such that $a_{i_p} = b_{j_p}$ for each $p \in [t]$.

\section{The compression operation}
\label{Compsection}
Compression operations have various useful properties. It is straightforward that for $i, j \in [n]$ and $\mathcal{A} \subseteq 2^{[n]}$,
\[|\Delta_{i,j}(\mathcal{A})| = |\mathcal{A}|.\] 
We will also need the following well-known basic result (see, for example, \cite[Lemma~2.1]{Borg11}).

\begin{lemma}\label{compcross} Let $\mathcal{A}$ and $\mathcal{B}$ be cross-$t$-intersecting subfamilies of $2^{[n]}$.\\
(i) For any $i, j \in [n]$, $\Delta_{i,j}(\mathcal{A})$ and $\Delta_{i,j}(\mathcal{B})$ are cross-$t$-intersecting subfamilies of $2^{[n]}$.\\
(ii) If $1 \leq t \leq r \leq s \leq n$, $\mathcal{A} \subseteq {[n] \choose \leq r}$, $\mathcal{B} \subseteq {[n] \choose \leq s}$, and $\mathcal{A}$ and $\mathcal{B}$ are compressed, then \[|A \cap B \cap [r+s-t]| \geq t\] for any $A \in \mathcal{A}$ and any $B \in \mathcal{B}$.
\end{lemma}
The only difference between Lemma~\ref{compcross} and \cite[Lemma~2.1]{Borg11} is that the latter is for $\mathcal{A} \subseteq {[n] \choose r}$ and $\mathcal{B} \subseteq {[n] \choose s}$; however, the former follows by the argument for the latter.

Suppose that a subfamily $\mathcal{A}$ of $2^{[n]}$ is not compressed. Then $\mathcal{A}$ can be transformed to a compressed family through left-compressions as follows. Since $\mathcal{A}$ is not compressed, we can find a left-compression that changes $\mathcal{A}$, and we apply it to $\mathcal{A}$ to obtain a new subfamily of $2^{[n]}$. We keep on repeating this (always applying a left-compression to the last family obtained) until we obtain a subfamily of $2^{[n]}$ that is invariant under any left-compression (such a point is indeed reached, because if $\Delta_{i,j}(\mathcal{F}) \neq \mathcal{F} \subseteq 2^{[n]}$ and $i < j$, then $0 < \sum_{G \in \Delta_{i,j}(\mathcal{F})} \sum_{b \in G} b < \sum_{F \in \mathcal{F}} \sum_{a \in F} a$).

Now consider $\mathcal{A}, \mathcal{B} \subseteq 2^{[n]}$ such that $\mathcal{A}$ and $\mathcal{B}$ are cross-$t$-intersecting. Then, by Lemma~\ref{compcross}, we can obtain $\mathcal{A}^*, \mathcal{B}^* \subseteq 2^{[n]}$ such that $\mathcal{A}^*$ and $\mathcal{B}^*$ are compressed and cross-$t$-intersecting, $|\mathcal{A}^*| = |\mathcal{A}|$, and $|\mathcal{B}^*| = |\mathcal{B}|$. Indeed, similarly to the above procedure, if we can find a left-compression that changes at least one of $\mathcal{A}$ and $\mathcal{B}$, then we apply it to both $\mathcal{A}$ and $\mathcal{B}$, and we keep on repeating this (always performing this on the last two families obtained) until we obtain $\mathcal{A}^*, \mathcal{B}^* \subseteq 2^{[n]}$ such that both $\mathcal{A}^*$ and $\mathcal{B}^*$ are invariant under any left-compression.

\section{Proof of the main result} \label{Weightedsection}

This section is dedicated to the proof of Theorem~\ref{xintweight}.

For the extremal cases of Theorem~\ref{xintweight}, we shall use the following two lemmas.

\begin{lemma}\label{complemma1} Let $1 \leq t \leq n$ and $T = [t]$. Let $\mathcal{H}$ be a compressed subfamily of $2^{[n]}$. Let $w \colon \mathcal{H} \rightarrow \mathbb{R}^+$ such that  $w(\delta_{i,j}(H)) \geq w(H)$ for every $H \in \mathcal{H}$ and every $i,j \in [n]$ with $i < j$. Then $w(\mathcal{H}(T')) \leq w(\mathcal{H}(T))$ for each $T' \in {[n] \choose t}$.
\end{lemma}
\textbf{Proof.}  Let $T' \in {[n] \choose t}$, and let $a_1, \dots, a_t$ be the elements of $T'$. Let $\mathcal{D}_0 = \mathcal{H}(T')$. Let $\mathcal{D}_1 = \Delta_{1,a_1}(\mathcal{D}_0), \dots, \mathcal{D}_t = \Delta_{t,a_t}(\mathcal{D}_{t-1})$. Since $\mathcal{H}$ is compressed, $\mathcal{D}_i \subseteq \mathcal{H}$ for each $i \in [t]$. It follows from the properties of $w$ and of left-compressions 
that $w(\mathcal{D}_0) \leq w(\mathcal{D}_1) \leq \dots \leq w(\mathcal{D}_t)$. Thus the result follows if we show that $\mathcal{D}_t \subseteq \mathcal{H}(T)$. 

Let $D_1 \in \mathcal{D}_1$. If $D_1 \notin \mathcal{D}_0$, then $D_1 = \delta_{1,a_1}(D) \neq D$ for some $D \in \mathcal{D}_0$, and hence $1 \in D_1$. Suppose $D_1 \in \mathcal{D}_0$, so $a_1 \in D_1$ by definition of $\mathcal{D}_0$. Since $D_1$ is also in $\mathcal{D}_1$, $\delta_{1,a_1}(D_1) \in \mathcal{D}_0$. Thus $a_1 \in \delta_{1,a_1}(D_1)$ by definition of $\mathcal{D}_0$. Since $a_1 \in D_1$, it follows that $1 \in D_1$.

Therefore, $1 \in H$ for each $H \in \mathcal{D}_1$, that is, $\mathcal{D}_1 \subseteq \mathcal{H}(\{1\})$. If $t = 1$, then we have $w(\mathcal{D}_0) \leq w(\mathcal{D}_1) \leq w(\mathcal{H}(\{1\})) = w(\mathcal{H}(T))$, as required.

Suppose $t \geq 2$. Since $\mathcal{D}_1 \subseteq \mathcal{H}(\{1\})$, we clearly have $1 \in H$ for each $H \in \mathcal{D}_2$. By an argument similar to that for $\mathcal{D}_1$, we also obtain that $2 \in H$ for each $H \in \mathcal{D}_2$. Continuing this way, we obtain that $1, \dots, t \in H$ for each $H \in \mathcal{D}_t$. Thus $\mathcal{D}_t \subseteq \mathcal{H}(T)$, as required.~\hfill{$\Box$}

\begin{lemma}\label{complemma2} Let $n$, $t$, $T$, $\mathcal{G}$, $\mathcal{H}$, $g$, and $h$ be as in Theorem~\ref{xintweight}. If $U \in \mathcal{A} \subseteq \mathcal{G}$, $V \in \mathcal{B} \subseteq  \mathcal{H}$, $|U| = |V| = t$, and $\mathcal{A}$ and $\mathcal{B}$ are cross-$t$-intersecting, then
\[g(\mathcal{A}) h(\mathcal{B}) \leq g(\mathcal{G}(T)) h(\mathcal{H}(T)),\]
and equality holds if and only if $\mathcal{A} = \mathcal{G}(V)$, $\mathcal{B} = \mathcal{H}(V)$, $g(\mathcal{G}(V)) = g(\mathcal{G}(T))$, and $h(\mathcal{H}(V)) = h(\mathcal{H}(T))$.  
\end{lemma}
\textbf{Proof.}  Since $\mathcal{A}$ and $\mathcal{B}$ are cross-$t$-intersecting, we have $U = V$, $\mathcal{A} \subseteq \mathcal{G}(V)$, and $\mathcal{B} \subseteq \mathcal{H}(V)$. By Lemma~\ref{complemma1}, $g(\mathcal{G}(V)) \leq g(\mathcal{G}(T))$ and $h(\mathcal{H}(V)) \leq h(\mathcal{H}(T))$. Hence the result.~\hfill{$\Box$} \\
\\
\textbf{Proof of Theorem~\ref{xintweight}.} We prove the result by induction on $n$. 

Consider the base case $n=t$. If $g(\mathcal{A})h(\mathcal{B}) \neq 0$, then $\mathcal{A} \neq \emptyset \neq \mathcal{B}$, and hence, since $\mathcal{A}$ and $\mathcal{B}$ are cross-$t$-intersecting, $\mathcal{A} = \{T\} = \mathcal{B}$.

Now consider $n \geq t+1$. Let $\mathcal{A} \subseteq \mathcal{G}$ and $\mathcal{B} \subseteq  \mathcal{H}$ such that $g(\mathcal{A}) h(\mathcal{B})$ is maximum under the condition that $\mathcal{A}$ and $\mathcal{B}$ are cross-$t$-intersecting. If $\mathcal{G}$ does not have a member of size at least $t$, then $\mathcal{A} = \emptyset$ or $\mathcal{B} = \emptyset$ (since $\mathcal{A}$ and $\mathcal{B}$ are cross-$t$-intersecting), and hence $g(\mathcal{A})h(\mathcal{B}) = 0 = g(\mathcal{G}(T))h(\mathcal{H}(T))$. Similarly, $g(\mathcal{A})h(\mathcal{B}) = 0 = g(\mathcal{G}(T))h(\mathcal{H}(T))$ if $\mathcal{H}$ does not have a member of size at least $t$. Therefore, we will assume that each of $\mathcal{G}$ and $\mathcal{H}$ has a member of size at least $t$. Since $\mathcal{G}$ and $\mathcal{H}$ are hereditary and compressed, we clearly have $T \in \mathcal{G}$ and $T \in \mathcal{H}$. Thus $g(\mathcal{G}(T)) > 0$ and $h(\mathcal{H}(T)) > 0$. Since $\mathcal{G}(T)$ and $\mathcal{H}(T)$ are cross-$t$-intersecting, it follows by the choice of $\mathcal{A}$ and $\mathcal{B}$ that 
\begin{equation} g(\mathcal{A}) h(\mathcal{B}) \geq g(\mathcal{G}(T)) h(\mathcal{H}(T)) > 0. \label{main0.1}
\end{equation} 
It follows that $\mathcal{A} \neq \emptyset \neq \mathcal{B}$. It also follows that no member of $\mathcal{A}$ is of size less than $t$,  because otherwise $\mathcal{B} = \emptyset$, contradicting (\ref{main0.1}). Similarly, no member of $\mathcal{B}$ is of size less than $t$.

As explained in Section~\ref{Compsection}, we apply left-compressions to $\mathcal{A}$ and $\mathcal{B}$ simultaneously until we obtain two compressed cross-$t$-intersecting families $\mathcal{A}^*$ and $\mathcal{B}^*$, respectively. Thus $|\mathcal{A}^*| = |\mathcal{A}|$ and $|\mathcal{B}^*| = |\mathcal{B}|$. Since $\mathcal{G}$ and $\mathcal{H}$ are compressed, $\mathcal{A}^* \subseteq \mathcal{G}$ and $\mathcal{B}^* \subseteq \mathcal{H}$. By (b), $g(\mathcal{A}) \leq g(\mathcal{A}^*)$ and $h(\mathcal{B}) \leq h(\mathcal{B}^*)$. By the choice of $\mathcal{A}$ and $\mathcal{B}$, we actually have $g(\mathcal{A}) = g(\mathcal{A}^*)$ and $h(\mathcal{B}) = h(\mathcal{B}^*)$. 

Suppose that $\mathcal{A}^* = \mathcal{G}(U)$ and $\mathcal{B}^* = \mathcal{H}(U)$ for some $U \in {[n] \choose t}$ such that $g(\mathcal{G}(U)) = g(\mathcal{G}(T))$ and $h(\mathcal{H}(U)) = h(\mathcal{H}(T))$. Then $g(\mathcal{G}(U)) > 0$ and $h(\mathcal{H}(U)) > 0$, so $\mathcal{G}(U) \neq \emptyset$  and $\mathcal{H}(U) \neq \emptyset$. Thus, since $\mathcal{G}$ and $\mathcal{H}$ are hereditary, $U \in \mathcal{A}^*$ and $U \in \mathcal{B}^*$. Hence $V \in \mathcal{A}$ for some $V \in {[n] \choose t}$, and $V' \in \mathcal{B}$ for some $V' \in {[n] \choose t}$. By Lemma~\ref{complemma2}, the result follows.

Therefore, we may assume that $\mathcal{A}$ and $\mathcal{B}$ are compressed.

We first consider $t+1 \leq n \leq t+2$. If $\mathcal{A}$ has a member of size $t$ and $\mathcal{B}$ has a member of size $t$, then the result follows by Lemma~\ref{complemma2}. Thus, without loss of generality, we may assume that no member of $\mathcal{A}$ is of size $t$.\medskip

Suppose $n = t+1$. Then $\mathcal{A} = \{[t+1]\} \subseteq \mathcal{G}(T) \backslash \{T\}$ (since $\mathcal{A} \neq \emptyset$ and $\mathcal{A} \cap {[n] \choose t} = \emptyset$) and $\mathcal{B} \subseteq \mathcal{H} \cap \left({[t+1] \choose t} \cup \{[t+1]\}\right)$. Thus we have
\begin{align} g(\mathcal{A}) h(\mathcal{B}) &\leq (g(\mathcal{G}(T)) - g(T))\left(h(\mathcal{H}(T)) + \sum_{H \in {[t+1] \choose t} \cap \mathcal{H} \backslash \{T\}} h(H)\right) \nonumber \\
&\leq (g(\mathcal{G}(T)) - g(T))\left(h(\mathcal{H}(T)) + th(T)\right) \quad \mbox{(by (b))} \nonumber \\
&= g(\mathcal{G}(T))h(\mathcal{H}(T)) + th(T)g(\mathcal{G}(T)) - g(T)\left(h(\mathcal{H}(T)) + th(T)\right) \nonumber \\
&\leq g(\mathcal{G}(T))h(\mathcal{H}(T)) + th(T)(g(T) + g([t+1])) - (t+1)g(T)h(T) \nonumber \\
&\leq g(\mathcal{G}(T))h(\mathcal{H}(T)) + th(T)\left(g(T) + \frac{g(T)}{t+u}\right) - (t+1)g(T)h(T) \quad \mbox{(by (a))}. \nonumber
\end{align}
Therefore, $g(\mathcal{A}) h(\mathcal{B}) \leq g(\mathcal{G}(T)) h(\mathcal{H}(T))$, and equality holds only if $u = 0$.\medskip

Suppose $n = t+2$. This case requires a number of observations followed by the separate treatment of a few sub-cases. 

Let $T_1 = [t+1]$, $T_1' = T \cup \{t+2\}$, and $T_2 = [t+2]$. For each $i \in \{t, t+1, t+2\}$ and each $\mathcal{F} \in \{\mathcal{A}, \mathcal{B}, \mathcal{G}, \mathcal{H}\}$, let $\mathcal{F}^{(i)} = \mathcal{F} \cap {[t+2] \choose i}$. Thus $\mathcal{A} = \mathcal{A}^{(t)} \cup \mathcal{A}^{(t+1)} \cup \mathcal{A}^{(t+2)}$ and $\mathcal{B} = \mathcal{B}^{(t)} \cup \mathcal{B}^{(t+1)} \cup \mathcal{B}^{(t+2)}$. Recall that $\mathcal{A}$ has no $t$-set, so $\mathcal{A}^{(t)} = \emptyset$. Since $\mathcal{A}^{(t+2)}, \mathcal{B}^{(t+2)} \subseteq {[t+2] \choose t+2} = \{[t+2]\}$, we have $\mathcal{A}^{(t+2)} \subseteq \mathcal{G}(T)$ and $\mathcal{B}^{(t+2)} \subseteq \mathcal{H}(T)$. Let 
\[\mathcal{A}_T = \mathcal{A} \cap \mathcal{G}(T), \quad \mathcal{A}_{\overline{T}} = \mathcal{A} \backslash \mathcal{A}_T, \quad \mathcal{B}_T = \mathcal{B} \cap \mathcal{H}(T), \quad \mathcal{B}_{\overline{T}} = \mathcal{B} \backslash \mathcal{B}_T.\] 
We have
\begin{equation} \mathcal{A}_T \subseteq \mathcal{G}(T) \backslash \{T\}, \quad \mathcal{A}_{\overline{T}} \subseteq \mathcal{G}^{(t+1)} \backslash \{T_1, T_1'\}, \quad \mathcal{B}_{\overline{T}} \subseteq \mathcal{H}^{(t)} \cup \mathcal{H}^{(t+1)} \backslash \{T_1, T_1'\}. \label{main0.2}
\end{equation}
Since $T \subsetneq T_1 \subsetneq T_2$, we have $g(T_1) \leq \frac{g(T)}{t+u}$, $g(T_2) \leq \frac{g(T_1)}{t+u} \leq \frac{g(T)}{(t+u)^2}$, $h(T_1) \leq \frac{h(T)}{t+u}$, and $h(T_2) \leq \frac{h(T_1)}{t+u} \leq \frac{h(T)}{(t+u)^2}$. 
Clearly, for each $U \in \mathcal{G}^{(t)}$, there is a composition of left-compressions that gives $T$ when applied to $U$, and hence $g(U) \leq g(T)$ by (b). Similarly, $h(V) \leq h(T)$ for each $V \in \mathcal{H}^{(t)}$, $g(U) \leq g(T_1)$ for each $U \in \mathcal{G}^{(t+1)}$, and $h(V) \leq h(T_1)$ for each $V \in \mathcal{H}^{(t+1)}$. 

Suppose $\mathcal{A}^{(t+1)} = \emptyset$. Then $\mathcal{A} = \{T_2\}$. Since $T_2 \in \mathcal{G}$ and $\mathcal{G}$ is hereditary, we have $T, T_1,T_1',T_2 \in \mathcal{G}(T)$. Thus 
\begin{align} g(\mathcal{G}(T)) &\geq g(T) + g(T_1) + g(T_1') + g(T_2) \nonumber \\
&\geq (t+u)^2g(T_2) + 2(t+u)g(T_2) + g(T_2) \nonumber \\
&\geq ((t+u)+1)^2g(T_2) =  (t+u+1)^2g(\mathcal{A}), \nonumber
\end{align} 
and hence $g(\mathcal{A}) \leq \frac{g(\mathcal{G}(T))}{(t+u+1)^2}$. Now
\begin{align} h(\mathcal{B}) &= h(\mathcal{B}_T) + h(\mathcal{B}_{\overline{T}}) \leq h(\mathcal{H}(T)) + \left( {t+2 \choose t} - 1 \right) h(T) + \left( {t+2 \choose t+1} - 2 \right)h(T_1) \nonumber \\
&\leq h(\mathcal{H}(T)) + \left( \frac{(t+2)(t+1)}{2} - 1 \right) h(T) + t\frac{h(T)}{t} = h(\mathcal{H}(T)) + \frac{t^2 + 3t + 2}{2}h(T) \nonumber \\
&\leq h(\mathcal{H}(T)) + \frac{t^2 + 3t + 2}{2}h(\mathcal{H}(T)) = \frac{t^2 + 3t + 4}{2}h(\mathcal{H}(T)). \nonumber
\end{align} 
Thus 
\begin{align} g(\mathcal{A})h(\mathcal{B}) &\leq \frac{t^2 + 3t + 4}{2(t+u+1)^2} g(\mathcal{G}(T)) h(\mathcal{H}(T)) \leq \frac{t^2 + 3t + 4}{2(t+1)^2} g(\mathcal{G}(T)) h(\mathcal{H}(T)). \nonumber
\end{align}
Hence $g(\mathcal{A})h(\mathcal{B}) \leq g(\mathcal{G}(T)) h(\mathcal{H}(T))$, and equality holds only if $u = 0$.

Suppose that $\mathcal{A}^{(t+1)}$ has at least $3$ sets. Let $U_1, U_2$, and $U_3$ be $3$ distinct sets in $\mathcal{A}^{(t+1)}$. Since $U_1, U_2, U_3 \in {[t+2] \choose t+1}$, no $t$-set is a subset of each of $U_1, U_2$, and $U_3$. Thus no $t$-set $t$-intersects each of $U_1, U_2$, and $U_3$, and hence $\mathcal{B}^{(t)} = \emptyset$. We have
\[ g(\mathcal{A}) \leq g(\mathcal{G}(T)) - g(T) + \left( {t+2 \choose t+1} - 2 \right)g(T_1) \leq g(\mathcal{G}(T)) - g(T) + t\frac{g(T)}{t+u} \leq g(\mathcal{G}(T)). \]
Similarly, $h(\mathcal{B}) \leq h(\mathcal{H}(T))$. Thus $g(\mathcal{A})h(\mathcal{B}) \leq g(\mathcal{G}(T)) h(\mathcal{H}(T))$, and equality holds only if $u = 0$.

We still need to consider $1 \leq |\mathcal{A}^{(t+1)}| \leq 2$, for which we need more detailed observations. Let $\mathcal{C}_0 = \mathcal{B}_{\overline{T}} \cap {[t+2] \choose t}$, $\mathcal{C}_1 = \mathcal{B}_{\overline{T}} \cap {[t+2] \choose t+1}$, $\mathcal{D}_0 = \mathcal{H}(T) \cap {[t+2] \choose t}$, and $\mathcal{D}_1 = \mathcal{H}(T) \cap {[t+2] \choose t+1}$. By (\ref{main0.2}), $\mathcal{B}_{\overline{T}} = \mathcal{C}_0 \cup \mathcal{C}_1$. If $\mathcal{H}^{(t+1)} \backslash \{T_1, T_1'\}$ has a set $V$, then $t+2 \in V$, and hence there is a composition of left-compressions that gives $T_1'$ when applied to $V$. Thus, if $\mathcal{H}^{(t+1)} \backslash \{T_1, T_1'\}$ is non-empty, then $T_1, T_1' \in \mathcal{H}(T)$ (as $\mathcal{H}$ is compressed, $T \subset T_1$, and $T \subset T_1'$), and hence we have
\begin{align} h(\mathcal{C}_1) &\leq \sum_{V \in \mathcal{H}^{(t+1)} \backslash \{T_1, T_1'\}} h(V) \quad \mbox{(by (\ref{main0.2}))} \nonumber \\
&\leq \sum_{V \in \mathcal{H}^{(t+1)} \backslash \{T_1, T_1'\}} h(T_1') \leq th(T_1') \leq t\frac{h(T_1) + h(T_1')}{2} = \frac{t}{2}|\mathcal{D}_1|. \nonumber
\end{align} 
If $\mathcal{H}^{(t+1)} \backslash \{T_1, T_1'\} = \emptyset$, then $\mathcal{C}_1 = \emptyset$, and hence we also have $h(\mathcal{C}_1) \leq \frac{t}{2}h(\mathcal{D}_1)$. With a slight abuse of notation, we set $g(T_1') = 0$ if $T_1' \notin \mathcal{G}$, and we set $g(T_2) = 0$  if $T_2 \notin \mathcal{G}$. Since $\mathcal{G}$ is hereditary, $T_1' \in \mathcal{G}$ if $T_2 \in \mathcal{G}$. Thus $g(T_1') \geq (t+u)g(T_2)$.

Suppose that $\mathcal{A}^{(t+1)}$ has exactly one set. Since $\mathcal{A}$ is compressed, $\mathcal{A}^{(t+1)} = \{T_1\}$. Thus $\mathcal{A} \subseteq \{T_1, T_2\}$, and hence $g(\mathcal{A}) \leq g(T_1) + g(T_2) = g(\mathcal{G}(T)) - g(T) - g(T_1')$. The $t$-sets that $t$-intersect $T_1$ are those in ${T_1 \choose t}$, so $\mathcal{B}^{(t)} \subseteq {T_1 \choose t}$, and hence
\[h(\mathcal{C}_0) \leq \sum_{V \in (\mathcal{H}^{(t)} \backslash \{T\}) \cap {T_1 \choose t}} h(V) \leq \left( {t+1 \choose t} - 1 \right) h(T) =  th(\mathcal{D}_0).\] 
We have
\begin{align} g(\mathcal{A})h(\mathcal{B}) &\leq (g(\mathcal{G}(T)) - g(T) - g(T_1'))(h(\mathcal{H}(T)) + h(\mathcal{B}_{\overline{T}})) \nonumber \\
&= g(\mathcal{G}(T))h(\mathcal{H}(T)) + g(\mathcal{G}(T)) h(\mathcal{B}_{\overline{T}}) - (g(T) + g(T_1'))(h(\mathcal{H}(T)) + h(\mathcal{B}_{\overline{T}})) \nonumber \\
&= g(\mathcal{G}(T))h(\mathcal{H}(T)) + (g(T_1) + g(T_2))h(\mathcal{B}_{\overline{T}}) - (g(T) + g(T_1'))h(\mathcal{H}(T)) \nonumber \\
&\leq g(\mathcal{G}(T))h(\mathcal{H}(T)) + \left(\frac{g(T)}{t+u} + \frac{g(T_1')}{t+u} \right)h(\mathcal{B}_{\overline{T}}) - (g(T) + g(T_1'))h(\mathcal{H}(T)) \nonumber \\
&= g(\mathcal{G}(T))h(\mathcal{H}(T)) + (g(T) + g(T_1')) \left( \frac{h(\mathcal{C}_0) + h(\mathcal{C}_1)}{t+u} - h(\mathcal{H}(T)) \right) \nonumber \\
&\leq g(\mathcal{G}(T))h(\mathcal{H}(T)) + (g(T) + g(T_1')) \left( \frac{th(\mathcal{D}_0) + \frac{t}{2}h(\mathcal{D}_1)}{t+u} - (h(\mathcal{D}_0) + h(\mathcal{D}_1)) \right). \nonumber 
\end{align}
Thus $g(\mathcal{A})h(\mathcal{B}) \leq g(\mathcal{G}(T))h(\mathcal{H}(T))$, and equality holds only if $u = 0$.

Suppose that $\mathcal{A}^{(t+1)}$ has exactly $2$ sets. Since $\mathcal{A}$ is compressed, $\mathcal{A}^{(t+1)} = \{T_1, T_1'\}$. The only $t$-set that $t$-intersects each of $T_1$ and $T_1'$ is $T$, so $\mathcal{B}^{(t)} \subseteq \{T\}$. Thus $\mathcal{C}_0 = \emptyset$, and hence $\mathcal{B}_{\overline{T}} = \mathcal{C}_1$. Since $\mathcal{D}_1 \subseteq \{T_1, T_1'\}$, $h(\mathcal{D}_1) \leq 2\frac{h(T)}{t+u} = 2\frac{h(\mathcal{D}_0)}{t+u}$. Since $h(\mathcal{H}(T)) \geq h(\mathcal{D}_0) + h(\mathcal{D}_1)$, $h(\mathcal{H}(T)) \geq \frac{t+u}{2}h(\mathcal{D}_1) + h(\mathcal{D}_1) = \left(\frac{t+u}{2} + 1\right)h(\mathcal{D}_1)$. We have
\begin{align} g(\mathcal{A})h(\mathcal{B}) &\leq (g(\mathcal{G}(T)) - g(T))(h(\mathcal{H}(T)) + h(\mathcal{C}_1)) \nonumber \\
&= g(\mathcal{G}(T))h(\mathcal{H}(T)) + g(\mathcal{G}(T)) h(\mathcal{C}_1) - g(T)(h(\mathcal{H}(T)) + h(\mathcal{C}_1)) \nonumber \\
&= g(\mathcal{G}(T))h(\mathcal{H}(T)) + (g(T_1) + g(T_1') + g(T_2))h(\mathcal{C}_1) - g(T)h(\mathcal{H}(T)) \nonumber \\
&\leq g(\mathcal{G}(T))h(\mathcal{H}(T)) + \left(\frac{2g(T)}{t+u} + \frac{g(T)}{(t+u)^2} \right)\frac{t}{2}h(\mathcal{D}_1) - g(T)\left(\frac{t+u}{2} + 1\right)h(\mathcal{D}_1) \nonumber \\
&= g(\mathcal{G}(T))h(\mathcal{H}(T)) + g(T)h(\mathcal{D}_1) \left( \frac{t}{t+u} + \frac{t}{2(t+u)^2} - \frac{t+u}{2} - 1 \right). \nonumber 
\end{align}
Thus $g(\mathcal{A})h(\mathcal{B}) \leq g(\mathcal{G}(T))h(\mathcal{H}(T))$, and equality holds only if $u = 0$.\medskip

Now consider $n \geq t+3$. 

Define $\mathcal{H}_0 = \{H \in \mathcal{H} \colon n \notin H\}$
and $\mathcal{H}_1 = \{H \backslash \{n\} \colon n \in H \in
\mathcal{H}\}$. Define $\mathcal{G}_0$, $\mathcal{G}_1$,
$\mathcal{A}_0$, $\mathcal{A}_1$, $\mathcal{B}_0$, and
$\mathcal{B}_1$ similarly. Since $\mathcal{A}$, $\mathcal{B}$,
$\mathcal{G}$, and $\mathcal{H}$ are compressed, we clearly have
that $\mathcal{A}_0$, $\mathcal{A}_1$, $\mathcal{B}_0$,
$\mathcal{B}_1$, $\mathcal{G}_0$, $\mathcal{G}_1$,
$\mathcal{H}_0$, and $\mathcal{H}_1$ are compressed. Since
$\mathcal{G}$ and $\mathcal{H}$ are hereditary, we clearly have
that $\mathcal{G}_0$, $\mathcal{G}_1$, $\mathcal{H}_0$, and
$\mathcal{H}_1$ are hereditary, $\mathcal{G}_1 \subseteq \mathcal{G}_0$, and $\mathcal{H}_1 \subseteq \mathcal{H}_0$. Obviously, we have $\mathcal{A}_0 \subseteq \mathcal{G}_0 \subseteq 2^{[n-1]}$,
$\mathcal{A}_1 \subseteq \mathcal{G}_1 \subseteq 2^{[n-1]}$,
$\mathcal{B}_0 \subseteq \mathcal{H}_0 \subseteq 2^{[n-1]}$, and
$\mathcal{B}_1 \subseteq \mathcal{H}_1 \subseteq 2^{[n-1]}$.

Let $h_0 : \mathcal{H}_0 \rightarrow \mathbb{R}^+$ such that $h_0(H) = h(H)$ for each $H \in \mathcal{H}_0$. Let $h_1 : \mathcal{H}_1 \rightarrow \mathbb{R}^+$ such that $h_1(H) = h(H \cup \{n\})$ for each $H \in \mathcal{H}_1$ (note that $H \cup \{n\} \in \mathcal{H}$ by definition of $\mathcal{H}_1$). By (a) and (b), we have the following consequences. For any $A, B \in \mathcal{H}_0$ with $A \subsetneq B$ and $|A| \geq t$,
\begin{equation} h_0(A) = h(A) \geq (t+u)h(B) = (t+u)h_0(B). \label{main1}
\end{equation}
For any $C \in \mathcal{H}_0$ and any $i,j \in [n-1]$ with $i < j$,
\begin{equation} h_0(\delta_{i,j}(C)) = h(\delta_{i,j}(C)) \geq h(C) = h_0(C). \label{main2}
\end{equation}
For any $A, B \in \mathcal{H}_1$ with $A \subsetneq B$ and $|A| \geq t$,
\begin{equation} h_1(A) = h(A \cup \{n\}) \geq (t+u)h(B \cup \{n\}) = (t+u)h_1(B). \label{main3}
\end{equation}
For any $C \in \mathcal{H}_1$ and any $i,j \in [n-1]$ with $i < j$,
\begin{equation} h_1(\delta_{i,j}(C)) = h(\delta_{i,j}(C) \cup \{n\}) = h(\delta_{i,j}(C \cup \{n\})) \geq h(C \cup \{n\}) = h_1(C). \label{main4}
\end{equation}
Therefore, we have shown that properties (a) and (b) are inherited by $h_0$ and $h_1$.

Since $\mathcal{B} = \mathcal{B}_0 \cup \mathcal{B}(\{n\})$, $\mathcal{B}_0 \cap \mathcal{B}(\{n\}) = \emptyset$, and $\mathcal{B}(\{n\}) = \{B \cup \{n\} \colon B \in \mathcal{B}_1\}$, we have
\begin{equation} h(\mathcal{B}) = {h}(\mathcal{B}_0) + h(\mathcal{B}(\{n\})) = {h_0}(\mathcal{B}_0) + {h_1}(\mathcal{B}_1). \label{main4.2}
\end{equation}
Along the same lines,
\begin{align} h(\mathcal{H}(T)) &= h(\mathcal{H}_0(T)) + h(\{H \in \mathcal{H} \colon T \cup \{n\} \subseteq H\}) \nonumber \\
&= {h_0}(\mathcal{H}_0(T)) + {h}(\{H \cup \{n\} \colon H \in \mathcal{H}_1(T)\}) \nonumber \\
&= {h_0}(\mathcal{H}_0(T)) + {h_1}(\mathcal{H}_1(T)). \label{main4.3}
\end{align}

Suppose $\mathcal{G}_1 = \emptyset$. Clearly, $\mathcal{A}$ and $\mathcal{B}_0$ are cross-$t$-intersecting. Since $\mathcal{G}_1 = \emptyset$, no set in $\mathcal{A}$ contains $n$, and hence $\mathcal{A}$ and $\mathcal{B}_1$ are cross-$t$-intersecting. Thus, by the induction hypothesis,
\begin{equation} g(\mathcal{A}){h_j}(\mathcal{B}_j) \leq g(\mathcal{G}(T)){h_j}(\mathcal{H}_j(T)) \quad \mbox{for each $j \in \{0, 1\}$}. \label{main4.4}
\end{equation}
Together with (\ref{main4.2}) and (\ref{main4.3}), this gives us
\begin{align} g(\mathcal{A}){h}(\mathcal{B})
&= g(\mathcal{A}){h_0}(\mathcal{B}_0) + g(\mathcal{A}){h_1}(\mathcal{B}_1) \nonumber \\
&\leq g(\mathcal{G}(T)){h_0}(\mathcal{H}_0(T)) + g(\mathcal{G}(T)){h_1}(\mathcal{H}_1(T)) \nonumber \\
&= g(\mathcal{G}(T)){h}(\mathcal{H}(T)). \nonumber 
\end{align}
By (\ref{main0.1}), equality holds throughout, and hence $g(\mathcal{A}){h}(\mathcal{B}) = g(\mathcal{G}(T)){h}(\mathcal{H}(T))$. Thus, in (\ref{main4.4}), we actually have equality. Suppose $u > 0$. By the induction hypothesis, for each $j \in \{0,1\}$, we have $\mathcal{A} = \mathcal{G}(V_j)$ and $\mathcal{B}_j = \mathcal{H}_j(V_j)$ for some $V_j \in {[n-1] \choose t}$ such that $g(\mathcal{G}(V_j)) = g(\mathcal{G}(T))$ and ${h_j}(\mathcal{H}_j(V_j)) = {h_j}(\mathcal{H}_j(T))$. Thus $g(\mathcal{G}(V_0)) > 0$, and hence $\mathcal{G}(V_0) \neq \emptyset$. Thus, since $\mathcal{G}$ is hereditary and $\mathcal{A} = \mathcal{G}(V_0)$, $V_0 \in \mathcal{A}$. Since $\mathcal{A}$ and $\mathcal{B}$ are cross-$t$-intersecting, $\mathcal{B} \subseteq \mathcal{H}(V_0)$. Since $\mathcal{A} = \mathcal{G}(V_0)$, and since $\mathcal{G}(V_0)$ and $\mathcal{H}(V_0)$ are cross-$t$-intersecting, it follows by the choice of $\mathcal{A}$ and $\mathcal{B}$ that $\mathcal{B} = \mathcal{H}(V_0)$. By (\ref{main0.1}), $\mathcal{H}(V_0) \neq \emptyset$. Since $\mathcal{H}$ is hereditary, $V_0 \in \mathcal{B}$. By Lemma~\ref{complemma2}, the result follows.\medskip

Now suppose that $\mathcal{G}_1$ is non-empty. If $\mathcal{H}_1 = \emptyset$, then the result follows by an argument similar to that for the case $\mathcal{G}_1 = \emptyset$ above. Thus we assume that $\mathcal{H}_1$ is non-empty. Since $\mathcal{G}_1 \subseteq \mathcal{G}_0$ and $\mathcal{H}_1 \subseteq \mathcal{H}_0$, $\mathcal{G}_0$ and $\mathcal{H}_0$ are non-empty too. 

Similarly to $h_0$ and $h_1$, let $g_0 : \mathcal{G}_0 \rightarrow \mathbb{R}^+$ such that $g_0(G) = g(G)$ for each $G \in \mathcal{G}_0$, and let $g_1 : \mathcal{G}_1 \rightarrow \mathbb{R}^+$ such that $g_1(G) = g(G \cup \{n\})$ for each $G \in \mathcal{G}_1$ (note that $G \cup \{n\} \in \mathcal{G}$ by definition of $\mathcal{G}_1$). Then properties (a) and (b) are inherited by $g_0$ and $g_1$ in the same way they are inherited by $h_0$ and $h_1$, as shown above; that is, similarly to (\ref{main1})--(\ref{main4}), we have the following. For any $A, B \in \mathcal{G}_0$ with $A \subsetneq B$ and $|A| \geq t$,
\begin{equation} g_0(A) \geq (t+u)g_0(B). \label{main5}
\end{equation}
For any $C \in \mathcal{G}_0$ and any $i,j \in [n-1]$ with $i < j$,
\begin{equation} g_0(\delta_{i,j}(C)) \geq g_0(C). \label{main6}
\end{equation}
For any $A, B \in \mathcal{G}_1$ with $A \subsetneq B$ and $|A| \geq t$,
\begin{equation} g_1(A) \geq (t+u)g_1(B). \label{main7}
\end{equation}
For any $C \in \mathcal{G}_1$ and any $i,j \in [n-1]$ with $i < j$,
\begin{equation} g_1(\delta_{i,j}(C)) \geq g_1(C). \label{main8}
\end{equation}

Similarly to (\ref{main4.2}) and (\ref{main4.3}), we have
\begin{gather} g(\mathcal{A}) = {g_0}(\mathcal{A}_0) + {g_1}(\mathcal{A}_1), \label{main4.1} \\
g(\mathcal{G}(T)) = {g_0}(\mathcal{G}_0(T)) + {g_1}(\mathcal{G}_1(T)). \label{main4.6}
\end{gather}

Clearly, $\mathcal{A}_0$ and $\mathcal{B}_0$ are cross-$t$-intersecting, $\mathcal{A}_0$ and $\mathcal{B}_1$ are cross-$t$-intersecting, and $\mathcal{A}_1$ and $\mathcal{B}_0$ are cross-$t$-intersecting.

Let us first assume that $\mathcal{A}_1$ and $\mathcal{B}_1$ are cross-$t$-intersecting too. Then, by the induction hypothesis,
\begin{equation} {g_i}(\mathcal{A}_i){h_j}(\mathcal{B}_j) \leq {g_i}(\mathcal{G}_i(T)){h_j}(\mathcal{H}_j(T)) \quad \mbox{for any $i, j \in \{0, 1\}$.} \label{main4.5}
\end{equation}
Together with (\ref{main4.2}), (\ref{main4.3}), (\ref{main4.1}), and (\ref{main4.6}), this gives us
\begin{align} g(\mathcal{A}){h}(\mathcal{B}) &= {g_0}(\mathcal{A}_0) {h_0}(\mathcal{B}_0) + {g_0}(\mathcal{A}_0) {h_1}(\mathcal{B}_1) + 
{g_1}(\mathcal{A}_1) {h_0}(\mathcal{B}_0) + {g_1}(\mathcal{A}_1){h_1}(\mathcal{B}_1) \nonumber \\
&\leq {g_0}(\mathcal{G}_0(T)){h_0}(\mathcal{H}_0(T)) + {g_0}(\mathcal{G}_0(T)){h_1}(\mathcal{H}_1(T)) + \nonumber \\
& \quad \mbox{ } {g_1}(\mathcal{G}_1(T)){h_0}(\mathcal{H}_0(T)) + {g_1}(\mathcal{G}_1(T)){h_1}(\mathcal{H}_1(T)) \nonumber \\
&= g(\mathcal{G}(T)){h}(\mathcal{H}(T)). \nonumber
\end{align}
By (\ref{main0.1}), equality holds throughout, and hence $g(\mathcal{A}){h}(\mathcal{B}) = g(\mathcal{G}(T)){h}(\mathcal{H}(T))$. Thus, in (\ref{main4.5}), we actually have equality. Suppose $u > 0$. By the induction hypothesis, we particularly have $\mathcal{A}_0 = \mathcal{G}_0(V_0)$ and $\mathcal{B}_0 = \mathcal{H}_0(V_0)$ for some $V_0 \in {[n-1] \choose t}$ such that ${g_0}(\mathcal{G}_0(V_0)) = {g_0}(\mathcal{G}_0(T))$ and ${h_0}(\mathcal{H}_0(V_0)) = {h_0}(\mathcal{H}_0(T))$. Recall that $T \in \mathcal{G}$, so $T \in \mathcal{G}_0$, and hence ${g_0}(\mathcal{G}_0(T)) > 0$. Thus ${g_0}(\mathcal{G}_0(V_0)) > 0$, and hence $\mathcal{G}_0(V_0) \neq \emptyset$. Since $\mathcal{G}_0$ is hereditary, it follows that $V_0 \in \mathcal{G}_0(V_0)$, and hence $V_0 \in \mathcal{A}$. Similarly, $V_0 \in \mathcal{B}$. By Lemma~\ref{complemma2}, the result follows.\medskip

We will now show that $\mathcal{A}_1$ and $\mathcal{B}_1$ are indeed cross-$t$-intersecting. Note that $\mathcal{A}_1$ and $\mathcal{B}_1$ are cross-$(t-1)$-intersecting. 

Suppose that $\mathcal{A}_1$ and $\mathcal{B}_1$ are not cross-$t$-intersecting. Then there exists $A^* \in \mathcal{A}_1$ such that $|A^* \cap B^*| = t-1$ for some $B^* \in \mathcal{B}_1$. Let $r = |A^*|+1$ and $s = n-r+t$. Let 
\begin{align} \mathcal{R} &= \{A \in \mathcal{A}_1 \colon |A| = r-1, \, |A \cap B| = t-1 \mbox{ for some } B \in \mathcal{B}_1\}, \nonumber \\ \mathcal{S} &= \{B \in \mathcal{B}_1 \colon |B| = s-1, \, |A \cap B| = t-1 \mbox{ for some } A \in \mathcal{A}_1\}. \nonumber
\end{align}
We have $A^* \in \mathcal{R}$. 

Consider any $R \in \mathcal{R}$ and $B \in \mathcal{B}_1$ such that $|R \cap B| < t$. Since $\mathcal{A}_1$ and $\mathcal{B}_1$ are cross-$(t-1)$-intersecting, $|R \cap B| = t-1$. We have
\begin{align} |B| &= |B \cap R| + |B \backslash R| = t-1 + |B \backslash R| \nonumber \\
&\leq t-1 + |[n-1] \backslash R| = t-1 + (n - 1) - (r-1) = s-1.\nonumber
\end{align}
Suppose $B \notin \mathcal{S}$. Then $|B| < s-1$. Thus we have 
\[|R \cup B| = |R| + |B| - |R \cap B| \leq r-1 + s-2 - t+1 = n-2,\] 
and hence $R \cup B \subsetneq [n-1]$. Let $c \in [n-1] \backslash (R \cup B)$. Since $B \in \mathcal{B}_1$, $B \cup \{n\} \in \mathcal{B}$. Let $C = \delta_{c,n}(B \cup \{n\})$. Since $c \notin B \cup \{n\}$, $C = B \cup \{c\}$. Since $\mathcal{B}$ is compressed, $C \in \mathcal{B}$. However, since $c \notin R \cup \{n\}$ and $|R \cap B| = t-1$, we have $|(R \cup \{n\}) \cap C| = t-1$, which is a contradiction as $\mathcal{A}$ and $\mathcal{B}$ are cross-$t$-intersecting, $R \cup \{n\} \in \mathcal{A}$, and $C \in \mathcal{B}$.

We have therefore shown that
\begin{equation}\mbox{for each $B \in \mathcal{B}_1$ such that $|R \cap B| < t$ for some $R \in \mathcal{R}$, $B \in \mathcal{S}$.} \label{main8.1}
\end{equation}
%
By a similar argument,
\begin{equation}\mbox{for each $A \in \mathcal{A}_1$ such that $|A \cap S| < t$ for some $S \in \mathcal{S}$, $A \in \mathcal{R}$.} \label{main8.2}
\end{equation}

For each $A \in \mathcal{A}_1 \cup \mathcal{B}_1$, let $A' = A \cup \{n\}$. Let $\mathcal{R}' = \{R' \colon R \in \mathcal{R}\}$ and $\mathcal{S}' = \{S' \colon S \in \mathcal{S}\}$. Since $\mathcal{R} \subseteq \mathcal{A}_1$ and $\mathcal{S} \subseteq \mathcal{B}_1$, $\mathcal{R}' \subseteq \mathcal{A}(\{n\})$ and $\mathcal{S}' \subseteq \mathcal{B}(\{n\})$. Let 
\[\mathcal{A}' = \mathcal{A} \cup \mathcal{R}, \quad \mathcal{A}'' = \mathcal{A} \backslash \mathcal{R}', \quad \mathcal{B}' = \mathcal{B} \backslash \mathcal{S}', \quad \mathcal{B}'' = \mathcal{B} \cup \mathcal{S}.\]
By (\ref{main8.1}), $\mathcal{A}'$ and $\mathcal{B}'$ are cross-$t$-intersecting. By (\ref{main8.2}), $\mathcal{A}''$ and $\mathcal{B}''$ are cross-$t$-intersecting. Since $\mathcal{G}$ and $\mathcal{H}$ are hereditary, and since $\mathcal{R}' \subseteq \mathcal{A} \subseteq \mathcal{G}$ and $\mathcal{S}' \subseteq \mathcal{B} \subseteq \mathcal{H}$, we have $\mathcal{R} \subseteq \mathcal{G}$ and $\mathcal{S} \subseteq \mathcal{H}$, and hence $\mathcal{A}', \mathcal{A}'' \subseteq \mathcal{G}$ and $\mathcal{B}', \mathcal{B}'' \subseteq \mathcal{H}$. 

Let $x = g(\mathcal{A})$, $x_1 = g(\mathcal{R}')$, $y = h(\mathcal{B})$, and $y_1 = h(\mathcal{S}')$. We use a double-counting argument to obtain $x \geq nx_1/r$ and $y \geq ny_1/s$. For any $R \in \mathcal{R}'$ and any set $A$ such that $A = \delta_{i,n}(R)$ for some $i \in [n] \backslash R$, we write $A < R$. If $A < R \in \mathcal{R}'$, then, since $\mathcal{A}$ is compressed and $n \in R \in \mathcal{A}$, we have $A \in \mathcal{A}_0$. For any $A \in \mathcal{A}_0$ and any $R \in \mathcal{R}'$, let 
\[ \chi(A,R) = \left\{ \begin{array}{ll}
1 & \mbox{if $A < R$};\\
0 & \mbox{otherwise.}
\end{array} \right. \]
Then $\sum_{A \in \mathcal{A}_0} \chi(A,R) = n-r$ for each $R \in \mathcal{R}'$. For each $A \in \mathcal{A}_0$, $\chi(A,R) = 1$ only if $|A| = |R|$ and $R = (A \backslash \{i\}) \cup \{n\}$ for some $i \in A$. Thus $\sum_{R \in \mathcal{R}'} \chi(A,R) \leq r$ for each $A \in \mathcal{A}_0$. We have 
\begin{align} (n-r)x_1 &= \sum_{R \in \mathcal{R}'} (n-r)g(R) = \sum_{R \in \mathcal{R}'} \sum_{A \in \mathcal{A}_0} \chi(A,R)g(R) = \sum_{A \in \mathcal{A}_0} \sum_{R \in \mathcal{R}'} \chi(A,R)g(R) \nonumber \\
& \leq \sum_{A \in \mathcal{A}_0} \sum_{R \in \mathcal{R}'} \chi(A,R)g(A) \quad \mbox{(by (b))} \nonumber \\
&\leq \sum_{A \in \mathcal{A}_0} rg(A) = rg(\mathcal{A}_0) = r(x - g(\mathcal{A}(\{n\}))) \leq r(x - x_1), \nonumber
\end{align}
so $x \geq nx_1/r$. Similarly, $y \geq ny_1/s$. 

Since $t-1 = |A^* \cap B^*| \leq |A^*| = r-1$, $r \geq t$. By (\ref{main8.1}), $B^* \in \mathcal{S}$. Since $t-1 = |A^* \cap B^*| \leq |B^*| = s-1$, $s \geq t$.

Suppose $r = t$. Then $s = n$. Thus $\mathcal{S}' = \{[n]\}$ and $\mathcal{S} = \{[n-1]\}$. Let $C^* = [t-1] \cup \{n\}$. Since $A^* \in \mathcal{A}_1$, $\mathcal{A}_1$ is compressed, and $|A^*| = r-1 = t-1$, we have $[t-1] \in \mathcal{A}_1$, and hence $C^* \in \mathcal{A}$. 

Suppose that there exists $D^* \in \mathcal{B}$ such that $D^* \neq [n]$. Since $\mathcal{A}$ and $\mathcal{B}$ are cross-$t$-intersecting, we have $|C^* \cap D^*| \geq t$. Thus $C^* \subseteq D^*$ as $|C^*| = t$. Since $D^* \neq [n]$, there exists $c \in [n]$ such that $c \notin D^*$. Thus $c \notin C^*$. Since $\mathcal{A}$ is compressed, $\delta_{c,n}(C^*) \in \mathcal{A}$. However, $|\delta_{c,n}(C^*) \cap D^*| = |C^* \backslash \{n\}| = t-1$, which is a contradiction as $\mathcal{A}$ and $\mathcal{B}$ are cross-$t$-intersecting.

Therefore, $\mathcal{B} = \{[n]\}$. Since $n-1 > t$, $h([n-1]) \geq (t+u)h([n]) \geq th([n])$. We have 
\begin{equation} h(\mathcal{B}'') = h([n]) + h([n-1]) \geq h([n]) + th([n]) = (t+1)h(\mathcal{B}) = (t+1)y. \nonumber 
\end{equation}
Since $x \geq nx_1/r = nx_1/t \geq (t+3)x_1/t$, $x_1 \leq tx/(t+3)$. We have 
\[g(\mathcal{A}'') = x-x_1 \geq x - \frac{tx}{t+3} = \frac{3x}{t+3}.\] 
Thus we obtain 
\[g(\mathcal{A}'')h(\mathcal{B}'') \geq \frac{3(t+1)xy}{t+3} > xy = g(\mathcal{A})h(\mathcal{B}),\] 
contradicting the choice of $\mathcal{A}$ and $\mathcal{B}$.

Therefore, $r \geq t+1$. Similarly, $s \geq t+1$. Since $r-1 \geq t$ and each set in $\mathcal{R}$ is of size $r-1$, $g(\mathcal{R}) \geq (t+u)g(\mathcal{R}')$. Similarly, $h(\mathcal{S}) \geq (t+u)g(\mathcal{S}')$. 

Consider any $R \in \mathcal{R}$. By definition of $\mathcal{R}$, there exists $B_R \in \mathcal{B}_1$ such that $|R \cap B_R| = t-1$. Thus $|R \cap {B_R}'| = t-1$. Since ${B_R}' \in \mathcal{B}$, and since $\mathcal{A}$ and $\mathcal{B}$ are cross-$t$-intersecting, $R \notin \mathcal{A}$. Therefore, $\mathcal{A} \cap \mathcal{R} = \emptyset$. Similarly, $\mathcal{B} \cap \mathcal{S} = \emptyset$. 

We have
\begin{align} &g(\mathcal{A}') = x + g(\mathcal{R}) \geq x + (t+u)g(\mathcal{R}') = x + (t+u)x_1,  \nonumber \\
&g(\mathcal{A}'') = x - x_1, \nonumber \\
&h(\mathcal{B}') = y - y_1, \nonumber \\
&{h}(\mathcal{B}'') = y + h(\mathcal{S}) \geq y + (t+u)h(\mathcal{S}') = y + (t+u)y_1.  \nonumber
\end{align}
By the choice of $\mathcal{A}$ and $\mathcal{B}$, 
\[g(\mathcal{A}'){h}(\mathcal{B}') \leq g(\mathcal{A}){h}(\mathcal{B}) \quad \mbox{and} \quad g(\mathcal{A}''){h}(\mathcal{B}'') \leq g(\mathcal{A}){h}(\mathcal{B}).\]
Thus we have
\begin{align} &(x + (t+u)x_1)(y - y_1) \leq xy \quad \mbox{and} \quad (x - x_1)(y + (t+u)y_1) \leq xy \nonumber \\
&\Rightarrow \; (t+u)x_1y \leq xy_1 + (t+u)x_1y_1 \quad \mbox{and} \quad (t+u)xy_1 \leq x_1y + (t+u)x_1y_1 \nonumber \\
&\Rightarrow \; (t+u)x_1y + (t+u)xy_1 \leq (xy_1 + (t+u)x_1y_1) + (x_1y + (t+u)x_1y_1) \nonumber \\
&\Rightarrow \; (t+u-1)(x_1y + xy_1) \leq 2(t+u)x_1y_1. \nonumber \\
&\Rightarrow \; (t+u-1) \left( x_1\frac{ny_1}{s} + \frac{nx_1}{r}y_1 \right) \leq 2(t+u)x_1y_1 \nonumber \\ 
& \Rightarrow \; (t+u-1)(r+s)x_1y_1n \leq 2(t+u)rsx_1y_1 \nonumber \\
&\Rightarrow \; (t+u-1)(n+t)n \leq 2(t+u)r(n-r+t). \nonumber
\end{align}
Using differentiation, we find that the maximum value of the function $f(z) = z(n-z+t)$ occurs at $z = \frac{n+t}{2}$. Thus $r(n-r+t) \leq \frac{n+t}{2}\left(n-\frac{n+t}{2}+t\right) = (n+t)^2/4$, and hence
\begin{align} &(t+u-1)(n+t)n \leq 2(t+u)(n+t)^2/4 \nonumber \\
&\Rightarrow \; 2(t+u-1)n \leq (t+u)(n+t) \nonumber \\
&\Rightarrow \; n \leq \frac{(t+u)t}{t+u-2}. \label{keyineq}
\end{align}
Since $u > \frac{6-t}{3}$, $\frac{(t+u)t}{t+u-2} < t+3$. Thus we have $n < t+3$, which is a contradiction. Hence the result.~\hfill{$\Box$}

\begin{rem}\label{proofremark} \emph{Note that the proof for the special case $n \leq t+2$ actually verifies the conjecture in Remark~\ref{remark2} for $n \leq t+2$. Also note that for $t \geq 3$, if we also settle the conjecture for $t + 3 \leq n \leq t+6$, then we can take $u=0$ and proceed for $n \geq t+7$ in exactly the same way we did for $n \geq t+3$, because again we obtain a contradiction to (\ref{keyineq}); thus, as mentioned in Remark~\ref{remark2}, this would settle the conjecture for $t \geq 3$.} 
\end{rem}

\section{Proof of Theorem~\ref{nrs}} \label{nrssection}

In this section, we use Theorem~\ref{xintweight} to prove Theorem~\ref{nrs}.

For a family $\mathcal{F}$ and an integer $r \geq 0$, we denote the families $\{F \in \mathcal{F} \colon |F| = r\}$ and $\{F \in \mathcal{F} \colon |F| \leq r\}$ by $\mathcal{F}^{(r)}$ and $\mathcal{F}^{(\leq r)}$, respectively.

We will need the following lemma only when dealing with the characterisation of the extremal structures in the proof of Theorem~\ref{nrs}. 

\begin{lemma}\label{charlemma} Let $t,r,s,u,m,$ and $n$ be as in Theorem~\ref{nrs}. Let $i,j \in [\max\{m,n\}]$ with $i < j$. Let $\mathcal{G} = 2^{[m]}$ and $\mathcal{H} = 2^{[n]}$. Let $\mathcal{A} \subseteq \mathcal{G}^{(r)}$ and $\mathcal{B} \subseteq \mathcal{H}^{(s)}$ such that $\mathcal{A}$ and $\mathcal{B}$ are cross-$t$-intersecting. Suppose that $\Delta_{i,j}(\mathcal{A}) = \mathcal{G}^{(r)}(T)$ and $\Delta_{i,j}(\mathcal{B}) = \mathcal{H}^{(s)}(T)$ for some $t$-element subset $T$ of $[\min\{m,n\}]$. Then $\mathcal{A} = \mathcal{G}^{(r)}(T')$ and $\mathcal{B} = \mathcal{H}^{(s)}(T')$ for some $t$-element subset $T'$ of $[\min\{m,n\}]$.
\end{lemma}
We prove the above lemma using the following special case of \cite[Lemma~5.6]{Borg9}.

\begin{lemma} \label{L comp} Let $t \geq 1$, $r \geq t+1$, $n \geq 2r-t+2$, and $i, j \in [n]$. Let $\mathcal{H} = 2^{[n]}$, and let $\mathcal{A}$ be a $t$-intersecting subfamily of $\mathcal{H}^{(r)}$. If $\Delta_{i,j}(\mathcal{A})$ is a largest $t$-star of $\mathcal{H}^{(r)}$, then $\mathcal{A}$ is a largest $t$-star of $\mathcal{H}^{(r)}$.
\end{lemma}
%
\textbf{Proof of Lemma~\ref{charlemma}.} We are given that $t \leq r \leq s$. 

Suppose $r=t$. Then $\Delta_{i,j}(\mathcal{A}) = \{T\}$. Thus $\mathcal{A} = \{T'\} = \mathcal{H}^{(r)}(T')$ for some $T' \in {[m] \choose t}$. Since $\mathcal{A}$ and $\mathcal{B}$ are cross-$t$-intersecting, $T' \subseteq B$ for all $B \in \mathcal{B}$. Thus $\mathcal{B} \subseteq \mathcal{H}^{(s)}(T')$. Since ${n-t \choose s-t} = |\mathcal{H}^{(s)}(T)| = |\Delta_{i,j}(\mathcal{B})| = |\mathcal{B}| \leq |\mathcal{H}^{(s)}(T')| = {n-t \choose s-t}$, $|\mathcal{B}| = {n-t \choose s-t}$. Hence $\mathcal{B} = \mathcal{H}^{(s)}(T')$. 

Now suppose $r \geq t+1$. Note that $T \backslash \{i\} \subseteq E$ for all $E \in \mathcal{A} \cup \mathcal{B}$.

Suppose that $\mathcal{A}$ is not $t$-intersecting. Then there exist $A_1, A_2 \in \mathcal{A}$ such that $|A_1 \cap A_2| \leq t-1$, and hence $T \nsubseteq A_l$ for some $l \in \{1,2\}$; we may assume that $l=1$. Thus, since $\Delta_{i,j}(\mathcal{A}) = \mathcal{H}^{(r)}(T)$, we have $A_1 \notin \Delta_{i,j}(\mathcal{A})$, $A_1 \neq \delta_{i,j}(A_1) \in \Delta_{i,j}(\mathcal{A})$, $\delta_{i,j}(A_1) \notin \mathcal{A}$ (because otherwise $A_1 \in \Delta_{i,j}(\mathcal{A})$), $i \in T$, $j \notin T$, $j \in A_1$, and $A_1 \cap T = T \backslash \{i\}$. Since $T \backslash \{i\} \subseteq A_1 \cap A_2$ and $|A_1 \cap A_2| \leq t-1$, $A_1 \cap A_2 = T \backslash \{i\}$. Thus $j \notin A_2$, and hence $A_2 = \delta_{i,j}(A_2)$. Since $\delta_{i,j}(A_2) \in \Delta_{i,j}(\mathcal{A}) = \mathcal{H}^{(r)}(T)$, $T \subseteq A_2$. Let $X = [n] \backslash (A_1 \cup A_2)$. We have
\begin{align} |X| &= n - |A_1 \cup A_2| = n - (|A_1| + |A_2| - |A_1 \cap A_2|) = n - 2r + t - 1 \nonumber \\
&\geq (t+u+2)(s-t)+r-1 - 2(r - t) - (t + 1) \geq (t+u)(s-t)-1 \nonumber \\
&> \left(2 + \frac{2t}{3} \right)(s-t)-1 = \left(1 + \frac{2t}{3} \right)(s-t) + s - (t+1).  \nonumber
\end{align} 
Since $t+1 \leq r \leq s$, we have $|X| > s-t$, and hence ${X \choose s-t} \neq \emptyset$. Let $C \in {X \choose s-t}$ and $D = C \cup T$. Then $D \in \mathcal{H}^{(s)}(T)$ and $D \cap A_1 = T \backslash \{i\}$, meaning that $D \in \Delta_{i,j}(\mathcal{B})$ and $|D \cap A_1| = t-1$. Since $\mathcal{A}$ and $\mathcal{B}$ are cross-$t$-intersecting, we obtain $D \notin \mathcal{B}$ and $(D \backslash \{i\}) \cup \{j\} \in \mathcal{B}$, which is a contradiction since $|((D \backslash \{i\}) \cup \{j\}) \cap A_2| = |T \backslash \{i\}| = t-1$.

Therefore, $\mathcal{A}$ is $t$-intersecting. Similarly, $\mathcal{B}$ is $t$-intersecting. Now $\mathcal{H}^{(r)}(T)$ is a largest $t$-star of $\mathcal{H}^{(r)}$, and $\mathcal{H}^{(s)}(T)$ is a largest $t$-star of $\mathcal{H}^{(s)}$. Since $t+1 \leq r \leq s$ and
\[\max\{m,n\} \geq (t+u+2)(s-t)+r-1 = (t+u)(s-t) + (2s - t + 2) + r - (t+1) - 2,\]
we have
\[\max\{m,n\} - (2s - t + 2) > \left(2 + \frac{2t}{3} \right)(s-t)-2 \geq \frac{2t}{3}(s-t) > 0,\]
and hence $\max\{m,n\} > 2s - t + 2$. By Lemma~\ref{L comp}, $\mathcal{A} = \mathcal{H}^{(r)}(T')$ for some $T' \in {[m] \choose t}$, and $\mathcal{B} = \mathcal{H}^{(s)}(T^*)$ for some $T^* \in {[n] \choose t}$. 

Suppose $T' \neq T^*$. Let $z \in T^* \backslash T'$. Since $m > 2s-t+2 > r$, we can choose $A' \in \mathcal{H}^{(r)}(T')$ such that $z \notin A'$. Since $n > 2s - t + 2 \geq r + s - t + 2 > r + s - t$ and $z \in T^* \backslash A'$, we can choose $B^* \in \mathcal{H}^{(s)}(T^*)$ such that $|A' \cap B^*| \leq t-1$; however, this is a contradiction since $\mathcal{A} = \mathcal{H}^{(r)}(T')$, $\mathcal{B} = \mathcal{H}^{(s)}(T^*)$, and $\mathcal{A}$ and $\mathcal{B}$ are cross-$t$-intersecting. Therefore, $T' = T^*$.~\hfill{$\Box$} \\
\\
%
\textbf{Proof of Theorem~\ref{nrs}.} If $\mathcal{A} = \emptyset$ or $\mathcal{B}= \emptyset$, then $|\mathcal{A}||\mathcal{B}| = 0$. Thus we assume that $\mathcal{A} \neq \emptyset$ and $\mathcal{B} \neq \emptyset$. Let $l = \max\{m,n\}$, so $\mathcal{A}, \mathcal{B} \subseteq 2^{[l]}$.

As explained in Section~\ref{Compsection}, we apply left-compressions to $\mathcal{A}$ and $\mathcal{B}$ simultaneously until we obtain two compressed cross-$t$-intersecting families $\mathcal{A}^*$ and $\mathcal{B}^*$, respectively. We have $\mathcal{A}^* \subseteq {[m] \choose r}$, $\mathcal{B}^* \subseteq {[n] \choose s}$, $|\mathcal{A}^*| = |\mathcal{A}|$, and $|\mathcal{B}^*| = |\mathcal{B}|$. In view of Lemma~\ref{charlemma}, we may therefore assume that $\mathcal{A}$ and $\mathcal{B}$ are compressed. Thus, by Lemma~\ref{compcross}(ii),
\begin{equation} \mbox{$|A \cap B \cap [r+s-t]| \geq t$ for any $A \in \mathcal{A}$ and any $B \in \mathcal{B}$.} \label{nrs1}
\end{equation}

Let $p = r+s-t$. Let $\mathcal{G} = {[p] \choose \leq r}$ and $\mathcal{H} = {[p] \choose \leq s}$. Let $g : \mathcal{G} \rightarrow \mathbb{N}$ such that $g(G) = {m-p \choose r-|G|}$ for each $G \in \mathcal{G}$. Let $h : \mathcal{H} \rightarrow \mathbb{N}$ such that $h(H) = {n-p \choose s-|H|}$ for each $H \in \mathcal{H}$. 

For every $F, G \in \mathcal{G}$ with $F \subsetneq G$ and $t \leq |F| = |G|-1$, we have
%
\begin{align} \frac{g(F)-(t+u)g(G)}{{m-p \choose r-|F|}} &= 1 - \frac{(t+u){m-p \choose r-|F|-1}}{{m-p \choose r-|F|}} = 1 - \frac{(t+u)(r-|F|)}{m-p-(r-|F|)+1} \nonumber \\
&= \frac{m-p-(t+u+1)(r-|F|)+1}{m-p-(r-|F|)+1} \nonumber \\
&\geq \frac{m-p-(t+u+1)(r-t)+1}{m-p-(r-|F|)+1}  \nonumber \\
&= \frac{m-(t+u+2)(r-t)-s+1}{m-p-(r-|F|)+1} \nonumber \\
&\geq \frac{(t+u+2)(s-t)+r-1-((t+u+2)(r-t)+s-1)}{m-p-(r-|F|)+1} \geq 0, \nonumber 
\end{align}
and hence $g(F) \geq (t+u)g(G)$. It follows that $g(F) \geq (t+u)g(G)$ for every $F, G \in \mathcal{G}$ with $F \subsetneq G$ and $|F| \geq t$. Similarly, $h(F) \geq (t+u)g(H)$ for every $F, H \in \mathcal{H}$ with $F \subsetneq H$ and $|F| \geq t$.

We have $g(\delta_{i,j}(G)) = g(G)$ for every $G \in \mathcal{G}$ and every $i,j \in [p]$. Similarly, $h(\delta_{i,j}(H)) = h(H)$ for every $H \in \mathcal{H}$ and every $i,j \in [p]$.

Let $\mathcal{C} = \{A \cap [p] \colon A \in \mathcal{A}\}$ and $\mathcal{D} = \{B \cap [p] \colon B \in \mathcal{B}\}$. Then $\mathcal{C} \subseteq \mathcal{G}$, $\mathcal{D} \subseteq \mathcal{H}$, and, by (\ref{nrs1}), $\mathcal{C}$ and $\mathcal{D}$ are cross-$t$-intersecting. Let $T = [t]$. By Theorem~\ref{xintweight},
\begin{equation} g(\mathcal{C}) h(\mathcal{D}) \leq g(\mathcal{G}(T)) h(\mathcal{H}(T)), \label{nrs2}
\end{equation}
and if $u > 0$, then equality holds only if $\mathcal{C} = \mathcal{G}(T')$ and $\mathcal{D} = \mathcal{H}(T')$ for some $T' \in {[p] \choose t}$.

We have 
\begin{align} &|\mathcal{A}| = \left| \bigcup_{i = 0}^r \left\{ A \in \mathcal{A} \colon A \cap [p] \in \mathcal{C}^{(i)} \right\} \right| \leq \sum_{i = 0}^r \left| \mathcal{C}^{(i)} \right| {m-p \choose r-i} = g(\mathcal{C}), \label{nrs3} \\
&|\mathcal{B}| = \left| \bigcup_{j = 0}^s \left\{ B \in \mathcal{B} \colon B \cap [p] \in \mathcal{D}^{(j)} \right\} \right| \leq \sum_{j = 0}^s \left| \mathcal{D}^{(j)} \right| {n-p \choose s-j} = h(\mathcal{D}), \label{nrs4}
\end{align}
and hence, by (\ref{nrs2}), 
\begin{equation} |\mathcal{A}||\mathcal{B}| \leq g(\mathcal{G}(T)) h(\mathcal{H}(T)). \label{nrs5}
\end{equation}
Now 
\begin{align} g(\mathcal{G}(T)) &= \sum_{i = t}^r \left| \mathcal{G}^{(i)}(T) \right| {m-p \choose r-i} = \left| \bigcup_{i = t}^r \left\{ A \in {[m] \choose r} \colon T \subseteq A, |A \cap [p]| = i \right\} \right|  \nonumber \\
&= \left| \left\{ A \in {[m] \choose r} \colon T \subseteq A \right\} \right| = {m-t \choose r-t} \nonumber
\end{align}
and, similarly, $h(\mathcal{H}(T)) = {n-t \choose s-t}$. Together with (\ref{nrs5}), this gives us
\[|\mathcal{A}||\mathcal{B}| \leq {m-t \choose r-t}{n-t \choose s-t},\]
as required.

Suppose $|\mathcal{A}||\mathcal{B}| = {m-t \choose r-t}{n-t \choose s-t}$ and $u > 0$. Then equality holds throughout in each of (\ref{nrs2})--(\ref{nrs5}), and hence $\mathcal{C} = \mathcal{G}(T')$ and $\mathcal{D} = \mathcal{H}(T')$ for some $T' \in {[p] \choose t}$. It follows that $\mathcal{A} \subseteq \left\{ A \in {[m] \choose r} \colon T' \subseteq A \right\}$ and $\mathcal{B} \subseteq \left\{ B \in {[n] \choose s} \colon T' \subseteq B \right\}$. Since $|\mathcal{A}||\mathcal{B}| = {m-t \choose r-t}{n-t \choose s-t}$, both inclusion relations are actually equalities, $T' \subseteq [m]$, and $T' \subseteq [n]$.~\hfill{$\Box$}

%
%
%
%

\section{Proof of Theorem~\ref{main}} \label{Proofmain}

In this section, we use Theorem~\ref{xintweight} to prove Theorem~\ref{main}.

We start by defining a compression operation for labeled sets. For any $x, y \in \mathbb{N}$, let
\[ \gamma_{x,y}(A) = \left\{ \begin{array}{ll}
(A \backslash \{(x,y)\}) \cup \{(x,1)\} & \mbox{if $(x,y) \in A$};\\
A & \mbox{otherwise}
\end{array} \right. \]
for any labeled set $A$, and let
\[\Gamma_{x,y}(\mathcal{A}) = \{\gamma_{x,y}(A) \colon A \in \mathcal{A}\} \cup \{A \in \mathcal{A} \colon \gamma_{x,y}(A) \in
\mathcal{A}\}\]
for any family $\mathcal{A}$ of labeled sets.

Note that $|\Gamma_{x,y}(\mathcal{A})| = |\mathcal{A}|$ and that if $\mathcal{A} \subseteq \mathcal{S}_{{\bf c},r}$, then $\Gamma_{x,y}(\mathcal{A}) \subseteq \mathcal{S}_{{\bf c},r}$. It is well known that if $\mathcal{A}$ and $\mathcal{B}$ are cross-$t$-intersecting families of labeled sets, then so are $\Gamma_{x,y}(\mathcal{A})$ and $\Gamma_{x,y}(\mathcal{B})$. We present a result that gives more than this.

For any IP sequence ${\bf c} = (c_1, \dots, c_n)$ and any $r \in [n]$, let $\mathcal{S}_{{\bf c},\leq r}$ denote the union $\bigcup_{i = 1}^r \mathcal{S}_{{\bf c},i}$.

\begin{lemma}\label{gamma} Let ${\bf c} = (c_1, \dots, c_m)$ and ${\bf d} = (d_1, \dots, d_n)$ be IP sequences. Let $x,y \in \mathbb{N}$, $y \geq 2$. Let $l = \max\{m,n\}$ and $h = \max\{c_m, d_n\}$. Let $V \subseteq [l] \times [2,h]$. Let $\mathcal{A} \subseteq \mathcal{S}_{{\bf c},\leq m}$ and $\mathcal{B} \subseteq \mathcal{S}_{{\bf d},\leq n}$ such that $|(A \cap B) \backslash V| \geq t$ for every $A \in \mathcal{A}$ and every $B \in \mathcal{B}$. Then $|(C \cap D) \backslash (V \cup \{(x,y)\})| \geq t$ for every $C \in \Gamma_{x,y}(\mathcal{A})$ and every $D \in \Gamma_{x,y}(\mathcal{B})$.
\end{lemma}
\textbf{Proof.}  Suppose $C \in \Gamma_{x,y}(\mathcal{A})$ and $D \in \Gamma_{x,y}(\mathcal{B})$. We first show that $|(C \cap D) \backslash
V| \geq t$. Let $C' = (C \backslash \{(x,1)\}) \cup \{(x,y)\}$ and $D' = (D \backslash \{(x,1)\}) \cup \{(x,y)\}$. If $C \in \mathcal{A}$ and $D \in \mathcal{B}$, then $|(C \cap D) \backslash V| \geq t$. If $C \notin \mathcal{A}$ and $D \notin \mathcal{B}$, then $(x,1) \in C \cap D$, $C' \in \mathcal{A}$, $D' \in \mathcal{B}$, and hence, since $(x,1) \notin V$, $|(C \cap D) \backslash V| \geq |(C' \cap D') \backslash V| \geq t$. Suppose $C \notin \mathcal{A}$ and $D \in \mathcal{B}$. Then $(x,1) \in C$ and $C' \in \mathcal{A}$. If $(x,y) \notin D$, then, since $C' \in \mathcal{A}$ and $D \in \mathcal{B}$, $t \leq |(C' \cap D) \backslash V| \leq |(C \cap D) \backslash V|$. If $(x,y) \in D$, then $\gamma_{x,y}(D) \in \mathcal{B}$ (because otherwise $D \notin \Gamma_{x,y}(\mathcal{B})$), and hence, since $C' \in \mathcal{A}$, $t \leq |(C' \cap \gamma_{x,y}(D)) \backslash V| = |(C \cap D) \backslash V|$. Similarly, if $C \in \mathcal{A}$ and $D \notin \mathcal{B}$, then $|(C \cap D) \backslash V| \geq t$.

Now suppose $(C \cap D) \backslash (V \cup \{(x,y)\}) < t$.
Since $|(C \cap D) \backslash V| \geq t$, $(x,y) \in C \cap D$. Thus $C, \gamma_{x,y}(C) \in \mathcal{A}$, $D, \gamma_{x,y}(D) \in \mathcal{B}$, and $|(C \cap \gamma_{x,y}(D)) \backslash V| = |(C \cap D) \backslash (V \cup \{(x,y)\})| < t$, a contradiction.~\hfill{$\Box$} 

\begin{cor}\label{deltacor} Let ${\bf c} = (c_1, \dots, c_m), {\bf d} = (d_1, \dots, d_n), l$, and $h$ be as in Lemma~\ref{gamma}. Let $\mathcal{A} \subseteq \mathcal{S}_{{\bf c},\leq m}$ and $\mathcal{B} \subseteq \mathcal{S}_{{\bf d},\leq n}$ such that $\mathcal{A}$ and $\mathcal{B}$ are cross-$t$-intersecting. Let
\begin{align} \mathcal{A}^* &= \Gamma_{l,h} \circ \dots \circ \Gamma_{l,2} \circ \dots \circ \Gamma_{2,h} \circ \dots \circ \Gamma_{2,2} \circ \Gamma_{1,h} \circ \dots \circ \Gamma_{1,2}(\mathcal{A}),\nonumber \\
\mathcal{B}^* &= \Gamma_{l,h} \circ \dots \circ \Gamma_{l,2} \circ \dots \circ \Gamma_{2,h} \circ \dots \circ \Gamma_{2,2} \circ \Gamma_{1,h} \circ \dots \circ \Gamma_{1,2}(\mathcal{B}). \nonumber
\end{align}
Then $|A \cap B \cap ([l] \times [1])| \geq t$ for every $A \in \mathcal{A}^*$ and every $B \in \mathcal{B}^*$.
\end{cor}
\textbf{Proof.}  Let $Z = [l] \times [2,h]$. By repeated application of Lemma~\ref{gamma},
$|(A \cap B) \backslash Z| \geq t$ for every $A \in \mathcal{A}^*$ and every $B \in \mathcal{B}^*$. The result follows since $(A \cap B) \backslash Z = A \cap B \cap ([l] \times [1])$.~\hfill{$\Box$} \\ 

The next lemma is needed for the characterisation of the extremal structures in Theorem~\ref{main}.

\begin{lemma}\label{ss_cross_int} Let ${\bf c} = (c_1, \dots, c_m), {\bf d} = (d_1, \dots, d_n), l$, and $h$ be as in Lemma~\ref{gamma}. Suppose $c_1 \geq 3$ and $d_1 \geq 3$. Let $r \in [m]$, $s \in [n]$, and $t \in [\min\{r,s\}]$. Let $\mathcal{A} \subseteq \mathcal{S}_{{\bf c},r}$ and $\mathcal{B} \subseteq \mathcal{S}_{{\bf d},s}$ such that $\mathcal{A}$ and $\mathcal{B}$ are cross-$t$-intersecting. Suppose $\Gamma_{x,y}(\mathcal{A}) = \mathcal{S}_{{\bf c},r}(T)$ and $\Gamma_{x,y}(\mathcal{B}) = \mathcal{S}_{{\bf d},s}(T)$ for some $(x,y) \in [l] \times [h]$ and some labeled set $T \in {[l] \times [h] \choose t}$. Then $\mathcal{A} = \mathcal{S}_{{\bf c},r}(T')$ and $\mathcal{B} = \mathcal{S}_{{\bf d},s}(T')$ for some labeled set $T' \in {[l] \times [h] \choose t}$.
\end{lemma}
\textbf{Proof.} The result is immediate if $\mathcal{A} = \Gamma_{x,y}(\mathcal{A})$ and $\mathcal{B} = \Gamma_{x,y}(\mathcal{B})$. Suppose $\mathcal{A} \neq \Gamma_{x,y}(\mathcal{A})$ or $\mathcal{B} \neq \Gamma_{x,y}(\mathcal{B})$. We may assume that $\mathcal{A} \neq \Gamma_{x,y}(\mathcal{A})$. Thus there exists $A_1 \in \mathcal{A} \backslash \Gamma_{x,y}(\mathcal{A})$ such that $\gamma_{x,y}(A_1) \in \Gamma_{x,y}(\mathcal{A}) \backslash \mathcal{A}$. Then $(x,1) \neq (x,y) \in A_1$ and $\gamma_{x,y}(A_1) = (A_1 \backslash \{(x,y)\}) \cup \{(x,1)\}$. 

Suppose $(x,1) \notin T$. Together with $\gamma_{x,y}(A_1) \in \Gamma_{x,y}(\mathcal{A}) = \mathcal{S}_{{\bf c},r}(T)$, this gives us $T \subseteq A_1$, which contradicts $A_1 \notin \Gamma_{x,y}(\mathcal{A})$.

Therefore, $(x,1) \in T$. Let $(a_1,b_1), \dots, (a_t,b_t)$ be the elements of $T$, where $(a_t,b_t) = (x,1)$. Let $T' = (T \backslash \{(x,1)\}) \cup \{(x,y)\}$. Since $\gamma_{x,y}(A_1) \in \mathcal{S}_{{\bf c},r}(T)$, we have $T' \subseteq A_1$, and hence $\mathcal{S}_{{\bf c},r}(T') \neq \emptyset$. Note that $|\mathcal{S}_{{\bf c},r}(T')| = |\mathcal{S}_{{\bf c},r}(T)|$.

Let $A^* \in \mathcal{S}_{{\bf c},r}(T')$. If $s > t$, then let $x_1, \dots, x_{s-t}$ be distinct elements of $[n] \backslash \{a_1, \dots, a_t\}$. For each $i \in [n]$, let $D_i = \{i\} \times [d_i]$. We are given that $3 \leq d_1 \leq \dots \leq d_n$. By definition of a labeled set, for each $i \in [n]$, we have $|A \cap D_i| \leq 1$ for all $A \in \mathcal{S}_{{\bf c},r}$. Thus $|D_i \backslash (A_1 \cup A^*)| \geq d_i - 2 \geq 1$ for each $i \in [n]$. If $s > t$, then let $(x_i,y_i) \in D_{x_i} \backslash (A_1 \cup A^*)$ for each $i \in [s-t]$, and let $B^* = T' \cup \{(x_1,y_1), \dots, (x_{s-t},y_{s-t})\}$. If $s = t$, then let $B^* = T'$. Thus $B^* \in \mathcal{S}_{{\bf d},s}(T')$. Since $\Gamma_{x,y}(\mathcal{B}) = \mathcal{S}_{{\bf d},s}(T)$, we have $B^* \in \mathcal{B}$ or $\gamma_{x,y}(B^*) \in \mathcal{B}$. However, $|\gamma_{x,y}(B^*) \cap A_1| = |T \cap A_1| = |T \backslash \{(x,1)\}| = t-1$, so $B^* \in \mathcal{B}$. Since $\Gamma_{x,y}(\mathcal{A}) = \mathcal{S}_{{\bf c},r}(T)$, we have $A^* \in \mathcal{A}$ or $\gamma_{x,y}(A^*) \in \mathcal{A}$. However, $|\gamma_{x,y}(A^*) \cap B^*| = t-1$, so $A^* \in \mathcal{A}$. 

We have therefore shown that $\mathcal{S}_{{\bf c},r}(T') \subseteq \mathcal{A}$. Since $|\mathcal{A}| = |\Gamma_{x,y}(\mathcal{A})| = |\mathcal{S}_{{\bf c},r}(T)| = |\mathcal{S}_{{\bf c},r}(T')|$, we actually have $\mathcal{A} = \mathcal{S}_{{\bf c},r}(T')$. Clearly, for each $L \in \mathcal{S}_{{\bf d},s}$ with $T' \nsubseteq L$, there exists $L' \in \mathcal{S}_{{\bf c},r}(T')$ such that $|L \cap L'| = |L \cap T'| < |T'| = t$. Thus, since $\mathcal{A} = \mathcal{S}_{{\bf c},r}(T')$, each set in $\mathcal{B}$ contains $T'$. Hence $\mathcal{B} \subseteq \mathcal{S}_{{\bf d},s}(T')$. Since $|\mathcal{B}| = |\Gamma_{x,y}(\mathcal{B})| = |\mathcal{S}_{{\bf d},s}(T)| = |\mathcal{S}_{{\bf d},s}(T')|$, we actually have $\mathcal{B} = \mathcal{S}_{{\bf d},s}(T')$.~\hfill{$\Box$} \\

The next lemma allows us to translate the setting in Theorem~\ref{main} to one given by Theorem~\ref{xintweight}.

\begin{lemma}\label{lemmaweighted} Let ${\bf c}$ be an IP sequence $(c_1, \dots, c_n)$. Let $r \in [n]$. Let $w \colon {[n] \choose \leq r} \rightarrow \mathbb{N}$ such that for each $A \in {[n] \choose \leq r}$,
\[w(A) = \left| \left\{ L \in \mathcal{S}_{{\bf c},r} \colon L \cap ([n] \times [1]) = A \times [1] \right\} \right|.\]
Then: \\
(i) $w(A) \geq (c_1-1)w(B)$ for every $A, B \in {[n] \choose \leq r}$ with $A \subsetneq B$. \\
(ii) $w(\delta_{i,j}(A)) \geq w(A)$ for every $A \in {[n] \choose \leq r}$ and every $i,j \in [n]$ with $i < j$.
\end{lemma}
\textbf{Proof.} (i) Let $A, B \in {[n] \choose \leq r}$ with $A \subsetneq B$. Let $B' = B \backslash A$. Thus $|B'| \geq 1$. For each $L \in \mathcal{S}_{{\bf c},r}$, let $\sigma(L) = \{x \in [n] \colon (x,y) \in L \mbox{ for some } y \in [c_i]\}$. We have
\begin{align} w(A) &\geq \left| \left\{ L \in \mathcal{S}_{{\bf c},r} \colon L \cap ([n] \times [1]) = A \times [1], \, B' \subseteq \sigma(L) \right\} \right| \nonumber \\
&= \sum_{E \in {[n] \backslash (A \cup B') \choose r - |A| - |B'|}} \prod_{b \in B'} (c_b - 1) \prod_{e \in E} (c_e - 1) \nonumber \\
&= \prod_{b \in B'} (c_b - 1) \left( \sum_{E \in {[n] \backslash B \choose r - |B|}} \prod_{e \in E} (c_e - 1) \right) \nonumber \\
&= w(B) \prod_{b \in B'} (c_b - 1) \geq (c_1-1)^{|B'|}w(B) \geq (c_1-1)w(B). \nonumber
\end{align}
(ii) Let $A \in {[n] \choose \leq r}$, and let $i,j \in [n]$ with $i < j$. Suppose $\delta_{i,j}(A) \neq A$. Then $j \in A$, $i \notin A$, and $\delta_{i,j}(A) = (A \backslash \{j\}) \cup \{i\}$. Let $B = A \backslash \{j\}$. Let
\begin{gather} \mathcal{E}_0 = {[n] \backslash (B \cup \{i, j\}) \choose r - |A|}, \nonumber \\
\mathcal{E}_1 = \left\{E \in {[n] \backslash (B \cup \{i\}) \choose r - |A|} \colon j \in E \right\}, \nonumber \\
\mathcal{E}_2 = \left \{E \in {[n] \backslash (B \cup \{j\}) \choose r - |A|} \colon i \in E \right\}. \nonumber
\end{gather}
We have
\begin{align} w(B \cup \{i\}) &= \sum_{E \in {[n] \backslash (B \cup \{i\}) \choose r - |A|}}\prod_{e \in E}(c_e - 1) \nonumber \\
&= \sum_{D \in \mathcal{E}_0}\prod_{d \in D}(c_d - 1) + \sum_{F \in \mathcal{E}_1}\prod_{f \in F}(c_f - 1) \nonumber \\
&\geq \sum_{D \in \mathcal{E}_0}\prod_{d \in D}(c_d - 1) + \sum_{F \in \mathcal{E}_1}\prod_{f \in F}(c_f - 1)\frac{c_i-1}{c_j-1} \quad \quad \mbox{(since $c_i \leq c_j$)} \nonumber \\
&= \sum_{D \in \mathcal{E}_0}\prod_{d \in D}(c_d - 1) + \sum_{F \in \mathcal{E}_2}\prod_{f \in F}(c_f - 1) \nonumber \\
&= \sum_{E \in {[n] \backslash (B \cup \{j\}) \choose r - |A|}}\prod_{e \in E}(c_e - 1) = w(B \cup \{j\}), \nonumber
\end{align}
and hence $w(\delta_{i,j}(A)) \geq w(A)$.~\hfill{$\Box$} \\ 
\\ 
\textbf{Proof of Theorem~\ref{main}.} Let $\mathcal{G} = {[m] \choose \leq r}$. Let $v \colon \mathcal{G} \rightarrow \mathbb{N}$ such that for each $G \in \mathcal{G}$,
\[v(G) = \left| \left\{ L \in \mathcal{S}_{{\bf c},r} \colon L \cap ([m] \times [1]) = G \times [1] \right\} \right|.\]
Let $\mathcal{H} = {[n] \choose \leq s}$. Let $w \colon \mathcal{H} \rightarrow \mathbb{N}$ such that for each $H \in \mathcal{H}$,
\[w(H) = \left| \left\{ L \in \mathcal{S}_{{\bf d},s} \colon L \cap ([n] \times [1]) = H \times [1] \right\} \right|.\]
Let $l = \max\{m,n\}$ and $h = \max\{c_m, d_n\}$. Let
\begin{align} \mathcal{A}^* &= \Gamma_{l,h} \circ \dots \circ \Gamma_{l,2} \circ \dots \circ \Gamma_{2,h} \circ \dots \circ \Gamma_{2,2} \circ \Gamma_{1,h} \circ \dots \circ \Gamma_{1,2}(\mathcal{A}),\nonumber \\
\mathcal{B}^* &= \Gamma_{l,h} \circ \dots \circ \Gamma_{l,2} \circ \dots \circ \Gamma_{2,h} \circ \dots \circ \Gamma_{2,2} \circ \Gamma_{1,h} \circ \dots \circ \Gamma_{1,2}(\mathcal{B}). \nonumber
\end{align}
Now let
\begin{align} \mathcal{C} &= \left\{G \in \mathcal{G} \colon E \cap ([m] \times [1]) = G \times [1] \mbox{ for some } E \in \mathcal{A}^* \right\}, \nonumber \\
\mathcal{D} &= \left\{H \in \mathcal{H} \colon F \cap ([n] \times [1]) = H \times [1] \mbox{ for some } F \in \mathcal{B}^* \right\}. \nonumber
\end{align}
Then $\mathcal{C} \subseteq \mathcal{G} \subseteq 2^{[l]}$, $\mathcal{D} \subseteq \mathcal{H} \subseteq 2^{[l]}$, and, by Corollary~\ref{deltacor}, $\mathcal{C}$ and $\mathcal{D}$ are cross-$t$-intersecting. We have 
\begin{align} \mathcal{A}^* &\subseteq \bigcup_{C \in \mathcal{C}} \{L \in \mathcal{S}_{{\bf c},r} \colon L \cap ([m] \times [1]) = C \times [1]\}, \label{maingen.01} \\
\mathcal{B}^* &\subseteq \bigcup_{D \in \mathcal{D}} \{L \in \mathcal{S}_{{\bf d},s} \colon L \cap ([n] \times [1]) = D \times [1]\}. \label{maingen.02}
\end{align}
Thus 
\begin{align} |\mathcal{A}^*| &\leq \sum_{C \in \mathcal{C}} v(C) = v(\mathcal{C}), \label{maingen.1} \\ 
|\mathcal{B}^*| &\leq \sum_{D \in \mathcal{D}} w(D) = w(\mathcal{D}). \label{maingen.12}
\end{align} 
Since $|\mathcal{A}| = |\mathcal{A}^*|$ and $|\mathcal{B}| = |\mathcal{B}^*|$, we therefore have 
\begin{align} |\mathcal{A}| &\leq v(\mathcal{C}), \label{maingen.15}\\ 
|\mathcal{B}| &\leq w(\mathcal{D}). \label{maingen.16}
\end{align}
Let $T_0 = [t]$. Let $\mathcal{I} = \mathcal{G}(T_0)$, $\mathcal{J} = \mathcal{H}(T_0)$, $\mathcal{X} = \mathcal{S}_{{\bf c},r}(T_0 \times [1])$, and $\mathcal{Y} = \mathcal{S}_{{\bf d},s}(T_0 \times [1])$. By Lemma~\ref{lemmaweighted} and Theorem~\ref{xintweight}, 
\begin{equation} v(\mathcal{C}) w(\mathcal{D}) \leq v(\mathcal{I}) w(\mathcal{J}). \label{maingen.2}
\end{equation}  
%
Now
\begin{align} v(\mathcal{I}) &= \left( \sum_{I \in \mathcal{I}} v(I) \right) = \left(\sum_{I \in \mathcal{I}} \left| \left\{L \in \mathcal{S}_{{\bf c},r} \colon L \cap ([m] \times [1]) = I \times [1] \right\} \right| \right) \nonumber \\
&= \left|\bigcup_{I \in \mathcal{I}} \left \{L \in \mathcal{S}_{{\bf c},r} \colon L \cap ([m] \times [1]) = I \times [1] \right\} \right| = |\mathcal{X}| \nonumber 
\end{align}
and, similarly, $w(\mathcal{J}) = |\mathcal{Y}|$. Together with (\ref{maingen.15})--(\ref{maingen.2}), this gives us $|\mathcal{A}||\mathcal{B}| \leq |\mathcal{X}||\mathcal{Y}|$, which establishes the first part of the theorem. 

We now prove the second part of the theorem. The sufficiency condition is trivial, so we prove the necessary condition.

Suppose $|\mathcal{A}||\mathcal{B}| = |\mathcal{X}||\mathcal{Y}|$ and $u > 0$. Then all the inequalities in (\ref{maingen.1})--(\ref{maingen.2}) are equalities. Having equality throughout in each of (\ref{maingen.1}) and (\ref{maingen.12}) implies that equality holds in each of (\ref{maingen.01}) and (\ref{maingen.02}). By Theorem~\ref{xintweight}, equality in (\ref{maingen.2}) gives us that $\mathcal{C} = \mathcal{G}(T_1)$ and $\mathcal{D} = \mathcal{H}(T_1)$ for some $T_1 \in {[l] \choose t}$. Together with equality in each of (\ref{maingen.01}) and (\ref{maingen.02}), this gives us that $\mathcal{A}^* = \mathcal{S}_{{\bf c},r}(T_2)$ and $\mathcal{B}^* = \mathcal{S}_{{\bf d},s}(T_2)$, where $T_2 = T_1 \times [1]$. By Lemma~\ref{ss_cross_int}, $\mathcal{A} = \mathcal{S}_{{\bf c},r}(T_3)$ and $\mathcal{B} = \mathcal{S}_{{\bf d},s}(T_3)$ for some $T_3 \in {[l] \times [h] \choose t}$. Since $|\mathcal{A}||\mathcal{B}| = |\mathcal{X}||\mathcal{Y}| > 0$, we clearly have $T_3 \in \mathcal{S}_{{\bf c},t} \cap \mathcal{S}_{{\bf d},t}$.~\hfill{$\Box$}

\section{Proof of Theorem~\ref{multithm}} \label{multisection}

In this section, we use Theorem~\ref{xintweight} to prove Theorem~\ref{multithm}. 

As in Section~\ref{nrssection}, for any family $\mathcal{F}$, $\mathcal{F}^{(r)}$ denotes $\{F \in \mathcal{F} \colon |F| = r\}$.  
For any $n,r \in \mathbb{N}$ and any family $\mathcal{A}$, let $M_{n,r,\mathcal{A}}$ denote the set $\{A \in M_{n,r} \colon {\rm S}_{A} \in \mathcal{A}\}$.

\begin{lemma}\label{multicomp} If $n, r \in \mathbb{N}$, $i,j \in [n]$, and $\mathcal{A} \subseteq 2^{[n]}$, then $|M_{n,r,\Delta_{i,j}(\mathcal{A})}| = |M_{n,r,\mathcal{A}}|$.
\end{lemma}
\textbf{Proof.} Let $\mathcal{B} = \Delta_{i,j}(\mathcal{A})$. Clearly, $|\mathcal{B}^{(p)}| = |\mathcal{A}^{(p)}|$ for each $p \in [n]$. We have
\begin{align} |M_{n,r,\mathcal{B}}| &= \sum_{B \in \mathcal{B}} |M_{n,r,\{B\}}| = \sum_{p = 1}^n \sum_{B \in \mathcal{B}^{(p)}} |M_{n,r,\{B\}}| = \sum_{p = 1}^n |\mathcal{B}^{(p)}| |M_{n,r,\{[p]\}}|\nonumber \\
&= \sum_{p = 1}^n |\mathcal{A}^{(p)}| |M_{n,r,\{[p]\}}| = \sum_{p = 1}^n \sum_{A \in \mathcal{A}^{(p)}} |M_{n,r,\{A\}}| = \sum_{A \in \mathcal{A}} |M_{n,r,\{A\}}| = |M_{n,r,\mathcal{A}}|, \nonumber
\end{align}
as required.~\hfill{$\Box$} \\ 
\\
\textbf{Proof of Theorem~\ref{multithm}.} 
Let $\mathcal{C} = \{{\rm S}_A \colon A \in \mathcal{A}\}$ and $\mathcal{D} = \{{\rm S}_B \colon B \in \mathcal{B}\}$. Clearly, $\mathcal{A} \subseteq M_{m,r,\mathcal{C}}$, $\mathcal{B} \subseteq M_{n,s,\mathcal{D}}$, and, since $\mathcal{A}$ and $\mathcal{B}$ are cross-$t$-intersecting, $M_{m,r,\mathcal{C}}$ and $M_{n,s,\mathcal{D}}$ are cross-$t$-intersecting. Thus we assume that 
\begin{equation}\mathcal{A} = M_{m,r,\mathcal{C}} \quad \mbox{and} \quad \mathcal{B} = M_{n,s,\mathcal{D}}. \label{mrs0}
\end{equation}
%

As explained in Section~\ref{Compsection}, we apply left-compressions to $\mathcal{C}$ and $\mathcal{D}$ simultaneously until we obtain two compressed cross-$t$-intersecting families $\mathcal{C}^*$ and $\mathcal{D}^*$, respectively. Since $\mathcal{C} \subseteq {[m] \choose \leq r}$ and $\mathcal{D} \subseteq {[n] \choose \leq s}$, we have $\mathcal{C}^* \subseteq {[m] \choose \leq r}$ and $\mathcal{D}^* \subseteq {[n] \choose \leq s}$. By Lemma~\ref{compcross}(ii),
\begin{equation} \mbox{$|C \cap D \cap [r+s-t]| \geq t$ for any $C \in \mathcal{C}^*$ and any $D \in \mathcal{D}^*$.} \label{mrs1}
\end{equation}

Let $p = r+s-t$. Let $\mathcal{G} = {[p] \choose \leq r}$ and $\mathcal{H} = {[p] \choose \leq s}$. Let $g : \mathcal{G} \rightarrow \mathbb{N}$ such that $g(G) = {m+r-p-1 \choose r-|G|}$ for each $G \in \mathcal{G}$. Let $h : \mathcal{H} \rightarrow \mathbb{N}$ such that $h(H) = {n+s-p-1 \choose s-|H|}$ for each $H \in \mathcal{H}$. 

For every $F, G \in \mathcal{G}$ with $F \subsetneq G$ and $t \leq |F| = |G|-1$, we have
\begin{align} \frac{g(F)-(t+u)g(G)}{{m+r-p-1 \choose r-|F|}} &= 1 - \frac{(t+u){m+r-p-1 \choose r-|F|-1}}{{m+r-p-1 \choose r-|F|}} = 1 - \frac{(t+u)(r-|F|)}{m-p+|F|} \nonumber \\
&= \frac{m-p+|F|-(t+u)(r-|F|)}{m-p+|F|} \nonumber \\
&\geq \frac{m-p+t-(t+u)(r-t)}{m-p+|F|} = \frac{m-(t+u+1)(r-t)-s+t}{m-p+|F|} \nonumber \\
&\geq \frac{(t+u+1)(s-t)+r-t-((t+u+1)(r-t)+s-t)}{m-p+|F|} \geq 0, \nonumber 
\end{align}
and hence $g(F) \geq (t+u)g(G)$. It follows that $g(F) \geq (t+u)g(G)$ for every $F, G \in \mathcal{G}$ with $F \subsetneq G$ and $|F| \geq t$. Similarly, $h(F) \geq (t+u)g(H)$ for every $F, H \in \mathcal{H}$ with $F \subsetneq H$ and $|F| \geq t$.

We have $g(\delta_{i,j}(G)) = g(G)$ for every $G \in \mathcal{G}$ and every $i,j \in [p]$. Similarly, $h(\delta_{i,j}(H)) = h(H)$ for every $H \in \mathcal{H}$ and every $i,j \in [p]$.

Let $\mathcal{E} = \{C \cap [p] \colon C \in \mathcal{C}^*\}$ and $\mathcal{F} = \{D \cap [p] \colon D \in \mathcal{D}^*\}$. Then $\mathcal{E} \subseteq \mathcal{G}$, $\mathcal{F} \subseteq \mathcal{H}$, and, by (\ref{mrs1}), $\mathcal{E}$ and $\mathcal{F}$ are cross-$t$-intersecting. Let $T = [t]$. By Theorem~\ref{xintweight},
\begin{equation} g(\mathcal{E}) h(\mathcal{F}) \leq g(\mathcal{G}(T)) h(\mathcal{H}(T)), \label{mrs2}
\end{equation}
and if $u > 0$, then equality holds only if $\mathcal{E} = \mathcal{G}(T')$ and $\mathcal{F} = \mathcal{H}(T')$ for some $T' \in {[p] \choose t}$.

By (\ref{mrs0}) and Lemma~\ref{multicomp},
\begin{align} |\mathcal{A}| &= |M_{m,r,\mathcal{C}^*}| \leq \left| \left\{A \in M_{m,r} \colon {\rm S}_A \cap [p] = E \mbox{ for some $E \in \mathcal{E}$} \right\} \right| \nonumber \\
&= \sum_{E \in \mathcal{E}} |\{A \in M_{m,r} \colon {\rm S}_A \cap [p] = E\}  \nonumber \\
&= \sum_{E \in \mathcal{E}} |\{(a_1, \dots, a_{r-|E|}) \colon a_1 \leq  \dots \leq a_{r-|E|}, \, a_1, \dots, a_{r-|E|} \in E \cup [p+1,m]\}|  \nonumber  \\
&= \sum_{E \in \mathcal{E}} |M_{|E|+m-p,r-|E|}| = \sum_{E \in \mathcal{E}} {m+r-p-1 \choose r-|E|} = g(\mathcal{E}).
\label{mrs3} 
\end{align}
Similarly, 
\begin{equation} |\mathcal{B}| \leq h(\mathcal{F}). \label{mrs4}
\end{equation}
By (\ref{mrs2})--(\ref{mrs4}), 
\begin{equation} |\mathcal{A}||\mathcal{B}| \leq g(\mathcal{G}(T)) h(\mathcal{H}(T)). \label{mrs5}
\end{equation}
Now, similarly to (\ref{mrs3}),
\begin{align} g(\mathcal{G}(T)) &= \left| \left\{A \in M_{n,r} \colon {\rm S}_A \cap [p] = E \mbox{ for some $E \in \mathcal{G}(T)$} \right\} \right| \nonumber \\
&= \left| \left\{A \in M_{n,r} \colon T \subseteq {\rm S}_A \right\} \right| = {m+r-t-1 \choose r-t}. \nonumber
\end{align}
Similarly, $h(\mathcal{H}(T)) = {n+s-t-1 \choose s-t}$. By (\ref{mrs5}), it follows that
\[|\mathcal{A}||\mathcal{B}| \leq {m+r-t-1 \choose r-t}{n+s-t-1 \choose s-t},\]
as required.

Suppose $|\mathcal{A}||\mathcal{B}| = {m+r-t-1 \choose r-t}{n+s-t-1 \choose s-t}$ and $u > 0$. Then equality holds throughout in each of (\ref{mrs2})--(\ref{mrs5}), and hence $\mathcal{E} = \mathcal{G}(T')$ and $\mathcal{F} = \mathcal{H}(T')$ for some $T' \in {[p] \choose t}$. Having equality throughout in (\ref{mrs3}) implies that $M_{m,r,\mathcal{C}^*} = \{A \in M_{m,r} \colon {\rm S}_A \cap [p] = E \mbox{ for some } E \in \mathcal{E}\} = \{A \in M_{m,r} \colon T' \subseteq {\rm S}_A\}$. Thus $T' \in \mathcal{C}^*$, and hence there exists $T_1 \in {[m] \choose t}$ such that $T_1 \in \mathcal{C}$. Similarly, there exists $T_2 \in {[n] \choose t}$ such that $T_2 \in \mathcal{D}$. Since $\mathcal{C}$ and $\mathcal{D}$ are cross-$t$-intersecting, we have $T_1 = T_2$, $\mathcal{C} \subseteq \{C \in {[m] \choose \leq r} \colon T_1 \subseteq C\}$, and $\mathcal{D} \subseteq \{D \in {[n] \choose \leq s} \colon T_1 \subseteq D\}$. Consequently, $\mathcal{A} \subseteq \left\{ A \in M_{m,r} \colon T_1 \subseteq {\rm S}_A \right\}$ and $\mathcal{B} \subseteq \left\{ B \in M_{n,s} \colon T_1 \subseteq {\rm S}_B \right\}$.  Since $|\mathcal{A}||\mathcal{B}| = {m+r-t-1 \choose r-t}{n+s-t-1 \choose s-t}$, both inclusion relations are actually equalities.~\hfill{$\Box$}

\section{The remaining cases} \label{conjecturesection}

Each of Theorems~\ref{nrs}, \ref{main}, and \ref{multithm} solves the particular cross-$t$-intersection problem under consideration for all cases where the ground sets are not smaller than a certain value dependent on $r$, $s$, and $t$. 
Solving any of these problems completely appears to be very difficult and would take this area of study to a significantly deeper level. 
We conjecture that the complete solutions are (\ref{conjecture1})--(\ref{conjecture3}) below.

For any $n \in \mathbb{N}$ and any $r, t, i, j \in \{0\} \cup [n]$ with $1 \leq t \leq r$ and $t+i+j \leq n$, let $\mathcal{M}_{n,r,t,i,j} = \{A \in {[n] \choose r} \colon |A \cap [t+i+j]| \geq t+i\}$. In \cite{F_t1}, Frankl conjectured that the size of a largest $t$-intersecting subfamily of ${[n] \choose r}$ is $\max \{ |\mathcal{M}_{m,r,t,i,i}| \colon i, j \in \{0\} \cup \mathbb{N}, \, t+2i \leq n\}$, and this was verified in \cite{AK1}. Hirschorn suggested an analogous conjecture \cite[Conjecture~4]{Hirschorn} for cross-$t$-intersecting families $\mathcal{A}$ and $\mathcal{B}$ with $\mathcal{A} \subseteq {[n] \choose r}$ and $\mathcal{B} \subseteq {[n] \choose s}$. Generalising Hirschorn's conjecture, we conjecture that if $m, n \in \mathbb{N}$, $r \in [m]$, $s \in [n]$, $t \in [\min\{r,s\}]$, $\mathcal{A} \subseteq {[m] \choose r}$, $\mathcal{B} \subseteq {[n] \choose s}$, and $\mathcal{A}$ and $\mathcal{B}$ are cross-$t$-intersecting, then
\begin{equation}|\mathcal{A}||\mathcal{B}| \leq \max \{ |\mathcal{M}_{m,r,t,i,j}| |\mathcal{M}_{n,s,t,j,i}| \colon i, j \in \{0\} \cup \mathbb{N}, \, t+i+j \leq \min\{m,n\} \}. \label{conjecture1}
\end{equation}

For any IP sequence ${\bf c} = (c_1, \dots, c_n)$ and any $r, t, i, j \in \{0\} \cup [n]$ with $1 \leq t \leq r$ and $t+i+j \leq n$, let $\mathcal{S}_{{\bf c},r,t,i,j} = \{A \in \mathcal{S}_{{\bf c},r} \colon |A \cap ([t+i+j] \times [1])| \geq t+i\}$. We conjecture that if ${\bf c} = (c_1, \dots, c_m)$ and ${\bf d} = (d_1, \dots, d_n)$ are IP sequences, $c_1 \geq 2$, $d_1 \geq 2$, $r \in [m]$, $s \in [n]$, $t \in [\min\{r,s\}]$, $\mathcal{A} \subseteq \mathcal{S}_{{\bf c},r}$, $\mathcal{B} \subseteq \mathcal{S}_{{\bf d},s}$, and $\mathcal{A}$ and $\mathcal{B}$ are cross-$t$-intersecting, then
\begin{equation} |\mathcal{A}||\mathcal{B}| \leq \max \{ |\mathcal{S}_{{\bf c},r,t,i,j}| |\mathcal{S}_{{\bf d},s,t,j,i}| \colon i, j \in \{0\} \cup \mathbb{N}, \, t+i+j \leq \min\{m,n\} \}. \label{conjecture2}
\end{equation}
This generalises \cite[Conjecture~3]{PT}, which is a conjecture for the case $r = s = m = n$.

For any $n \in \mathbb{N}$ and any $r, t, i, j \in \{0\} \cup [n]$ with $1 \leq t \leq r$ and $t+i+j \leq n$, let $M_{n,r,t,i,j} = \{A \in M_{n,r} \colon |{\rm S}_A \cap [t+i+j]| \geq t+i\}$. We conjecture that if $m, n, r, s \in \mathbb{N}$, $t \in [\min\{r,s\}]$, $\mathcal{A} \subseteq M_{m,r}$, $\mathcal{B} \subseteq M_{n,s}$, and $\mathcal{A}$ and $\mathcal{B}$ are cross-$t$-intersecting, then
\begin{equation}|\mathcal{A}||\mathcal{B}| \leq \max \{ |M_{m,r,t,i,j}| |M_{n,s,t,j,i}| \colon i, j \in \{0\} \cup \mathbb{N}, \, t+i+j \leq \min\{m,n\} \}. \label{conjecture3}
\end{equation}


\begin{thebibliography}{}

\bibitem{AK1} R. Ahlswede and L.H. Khachatrian, The complete
intersection theorem for systems of finite sets, European J.
Combin. 18 (1997), 125--136.

\bibitem{AK2} R. Ahlswede and L.H. Khachatrian, The diametric theorem
in Hamming spaces---Optimal anticodes, Adv. Appl. Math. 20
(1998), 429--449.

\bibitem{AC} M.O. Albertson and K.L. Collins, Homomorphisms of $3$-chromatic graphs, Discrete Math. 54 (1985), 127--132.



\bibitem{Bey3} C. Bey, On cross-intersecting families of sets, Graphs Combin. 21 (2005), 161--168. 

\bibitem{Bey1} C. Bey, The Erd\H{o}s--Ko--Rado bound for the function lattice, Discrete Appl. Math. 95 (1999), 115--125.



\bibitem{BorgBLMS} P. Borg, A cross-intersection theorem for subsets of a set, Bull. London. Math. Soc. 47 (2015), 248--256.

\bibitem{Borg4} P. Borg, A short proof of a cross-intersection
theorem of Hilton, Discrete Math. 309 (2009), 4750--4753.

\bibitem{Borg3} P. Borg, Cross-intersecting families of
permutations, J. Combin. Theory Ser. A 117 (2010), 483--487.

\bibitem{Borg2} P. Borg, Cross-intersecting families of partial
permutations, SIAM J. Disc. Math. 24 (2010), 600--608.


\bibitem{Borg5} P. Borg, Cross-intersecting sub-families of
hereditary families, J. Combin. Theory Ser. A 119 (2012), 871--881.

\bibitem{Borg9} P. Borg, Extremal $t$-intersecting sub-families of
hereditary families, J. London Math. Soc. 79 (2009), 167--185.

\bibitem{Borg} P. Borg, Intersecting and cross-intersecting families of labeled sets, Electron. J. Combin. 15 (2008), note paper N9.

\bibitem{Borg7} P. Borg, Intersecting families of sets and
permutations: a survey, in: Advances in Mathematics Research
(A.R. Baswell Ed.), Volume 16, Nova Science Publishers, Inc.,
2011, pp. 283--299.



\bibitem{Borg11} P. Borg, The maximum product of sizes of cross-$t$-intersecting uniform families, Australas. J. Combin. 60 (2014), 69--78.

\bibitem{Borg8} P. Borg, The maximum sum and the maximum product of sizes of cross-intersecting families, European J. Combin. 35 (2014), 117--130.

\bibitem{BL2} P. Borg and I. Leader, Multiple cross-intersecting
families of signed sets, J. Combin. Theory Ser. A 117 (2010),
583--588.

\bibitem{CK} P.J. Cameron and C.Y. Ku, Intersecting families of permutations, European J. Combin. 24 (2003), 881--890.

\bibitem{D} D.E. Daykin, Erd\H os--Ko--Rado from Kruskal--Katona,
J. Combin. Theory Ser. A 17 (1974), 254--255.

\bibitem{DF} M. Deza and P. Frankl, The Erd\H os--Ko--Rado theorem---22 years later, SIAM J. Algebraic Discrete Methods 4 (1983), pp. 419--431.



\bibitem{EKR} P. Erd\H os, C. Ko and R. Rado, Intersection
theorems for systems of finite sets, Quart. J. Math. Oxford (2)
12 (1961), 313--320.

\bibitem{F2} P. Frankl, Extremal set systems, in: R.L. Graham, M. Gr\"{o}tschel and L. Lov\'{a}sz (Eds.), Handbook of Combinatorics, Vol. 2, Elsevier, Amsterdam, 1995, pp. 1293--1329.

\bibitem{F_t1} P. Frankl, The Erd\H os--Ko--Rado Theorem is true
for $n = ckt$, Proc. Fifth Hung. Comb. Coll., North-Holland,
Amsterdam, 1978, pp. 365--375.

\bibitem{F} P. Frankl, The shifting technique in extremal set theory, in: C. Whitehead (Ed.), Surveys in Combinatorics, Cambridge Univ. Press, London/New York, 1987, pp. 81--110.


\bibitem{FT2} P. Frankl and N. Tokushige, The Erd\H os--Ko--Rado theorem for integer sequences, Combinatorica 19 (1999), 55--63.

\bibitem{FLST} P. Frankl, S.J. Lee, M. Siggers and N. Tokushige, An Erd\H{o}s--Ko--Rado theorem for cross $t$-intersecting families, J. Combin. Theory Ser. A 128 (2014), 207--249.

\bibitem{FGV} Z. F\"{u}redi, D. Gerbner and  M. Vizer, A discrete isodiametric result: the Erd\H{o}s--Ko--Rado theorem for multisets, European  J. Combin. 48 (2015), 224--233.


\bibitem{H} A.J.W. Hilton, An intersection theorem for a collection of families of subsets of a finite set, J. London Math. Soc. (2) 15 (1977), 369--376.

\bibitem{HM} A.J.W. Hilton and E.C. Milner, Some intersection theorems for systems of finite sets, Quart. J. Math. Oxford (2) 18 (1967), 369--384.

\bibitem{Hirschorn} J. Hirschorn, Asymptotic upper bounds on the shades of $t$-intersecting families, arXiv:0808.1434 [math.CO].

\bibitem{HST} F.C. Holroyd, C. Spencer and J. Talbot, Compression and Erd\H os--Ko--Rado graphs, Discrete Math. 293 (2005), 155--164.

\bibitem{HT} F.C. Holroyd and J. Talbot, Graphs with the Erd\H os--Ko--Rado property, Discrete Math. 293 (2005), 165--176.

\bibitem{K} G.O.H. Katona, A simple proof of the Erd\H os--Chao
Ko--Rado theorem, J. Combin. Theory Ser. B 13 (1972), 183--184.

\bibitem{Ka} G.O.H. Katona, A theorem of finite sets, in:
Theory of Graphs, Proc. Colloq. Tihany, Akad\'{e}miai Kiad\'{o}, 1968,
pp. 187--207.

\bibitem{Kat} G.O.H. Katona, Intersection theorems for systems of finite sets, Acta Math. Acad. Sci. Hungar. 15 (1964), 329--337.


\bibitem{Kr} J.B. Kruskal, The number of simplices in a
complex, in: Mathematical Optimization Techniques, University of
California Press, Berkeley, California, 1963, pp. 251--278.


\bibitem{MT2} M. Matsumoto and N. Tokushige, A generalization of the Katona theorem for cross $t$-intersecting families, Graphs Combin. 5 (1989), 159--171.

\bibitem{MT} M. Matsumoto and N. Tokushige, The exact bound in the
Erd\H os--Ko--Rado theorem for cross-intersecting families, J.
Combin. Theory Ser. A 52 (1989), 90--97.

\bibitem{MP} K. Meagher and A. Purdy, An Erd\H os--Ko--Rado theorem
for multisets, Electron. J. Combin. 18(1) (2011), paper P220.


\bibitem{Moon} A. Moon, An analogue of the Erd\H os--Ko--Rado theorem
for the Hamming schemes $H(n,q)$, J. Combin. Theory Ser. A 32
(1982), 386--390.

\bibitem{PT} J. Pach and G. Tardos, Cross-intersecting families of vectors, Graphs Combin. 31 (2015), 477--495.

\bibitem{Pyber} L. Pyber, A new generalization of the Erd\H os--Ko--Rado theorem, J. Combin. Theory Ser. A 43 (1986), 85--90.

\bibitem{ST} S. Suda and H. Tanaka, A cross-intersection theorem for  vector spaces based on semidefinite programming, Bull. London. Math. Soc. 46 (2014), 342--348.

\bibitem{T} J. Talbot, Intersecting families of separated sets, J. London Math. Soc. 68 (2003), 37--51.

\bibitem{Tok3} N. Tokushige, Cross $t$-intersecting integer sequences from weighted Erd\H{o}s--Ko--Rado, Combin. Prob. Comput. 22 (2013), 622--637. 

\bibitem{Tok1} N. Tokushige, On cross $t$-intersecting families of
sets, J. Combin. Theory Ser. A 117 (2010), 1167--1177.

\bibitem{Tok2} N. Tokushige, The eigenvalue method for cross $t$-intersecting families, J. Algebr. Comb. 38 (2013), 653--662.

\bibitem{WZ} J. Wang and H. Zhang, Cross-intersecting families and primitivity of symmetric systems, J. Combin. Theory Ser. A 118  (2011), 455--462.

\bibitem{W} R.M. Wilson, The exact bound in the Erd\H os--Ko--Rado
theorem, Combinatorica 4 (1984), 247--257.

\bibitem{Zhang} H. Zhang, Cross-intersecting families of labeled sets, Electron. J. Combin. 20(1) (2013), paper P17.


\end{thebibliography}
\end{document}